\documentclass[a4paper,11pt]{article}
\usepackage[latin1]{inputenc}
\usepackage[T1]{fontenc}
\usepackage[english]{babel}
\usepackage{layout}
\usepackage{amsmath}
\usepackage{amsfonts}
\usepackage{amssymb}
\usepackage{graphicx}
\usepackage{float}
\usepackage{multicol}
\usepackage{color}
\definecolor{gris25}{gray}{0.55}
\usepackage{shadow}
\usepackage{natbib}
\bibpunct{(}{)}{;}{a}{,}{,}

\newcommand{\be}{\begin{equation}}
\newcommand{\ee}{\end{equation}}
\newcommand{\bd}{\begin{displaymath}}
\newcommand{\ed}{\end{displaymath}}

\newcommand{\ba}{\begin{eqnarray}}
\newcommand{\ea}{\end{eqnarray}}
\newcommand{\ban}{\begin{eqnarray*}}
\newcommand{\ean}{\end{eqnarray*}}

\newcommand{\R} {\mathbb{R}}
\newcommand{\PP} {\mathbb{P}}
\newcommand{\E} {\mathbb{E}}
\newcommand{\N} {\mathbb{N}}

\newcommand{\cov} {\mbox{cov}}
\newcommand{\var} {\mbox{Var}}
\newcommand{\ds}{\displaystyle}

\DeclareMathOperator*{\argmax}{arg\,max}

\newcommand{\Keywords}[1]{\par\noindent
{\small{\textbf{Keywords}\/}: #1}}

\newcommand{\Abstract}[1]{\par\noindent
{\small{\textbf{Abstract}\/}: #1}}

\newcommand{\Acknowledgments}[1]{\par\noindent
{\large{\textbf{Acknowledgments}\/}. #1}}

\usepackage{theorem}
\theoremstyle{break}
\newtheorem{theo}{Theorem}
\newtheorem{lemm}{Lemma}
\newtheorem{coro}{Corollary}
\newtheorem{Remark}{Remark}
\theorembodyfont{\rmfamily}

\numberwithin{equation}{section}

\numberwithin{figure}{section}

\numberwithin{table}{section}

\numberwithin{theo}{section}

\numberwithin{lemm}{section}

\numberwithin{coro}{section}

\numberwithin{Remark}{section}

\begin{document}
\title{Off-line detection of multiple change points by the Filtered Derivative with p-Value method}

\date{ }
\maketitle \vspace{-1.5cm}
 \begin{center}
 Pierre, R.~BERTRAND${}^{1,2}$, Mehdi FHIMA${}^{2}$ and Arnaud GUILLIN${}^{2}$\\
 ${}^{1}$  {\it INRIA Saclay} \\
${}^{2}$ {\it Laboratoire de Math\'ematiques, UMR CNRS 6620\\
\& Universit\'e de Clermont-Ferrand II, France}
\end{center}

\Abstract{This paper deals with off-line detection of  change
points, for time series of independent observations, when the
number of change points is  unknown. We propose a sequential
analysis method with linear time and memory complexity. Our
method is based at first step, on Filtered Derivative method
which detects the right change points as well as the false ones. We
improve the Filtered Derivative method by adding a second step in
which we  compute the p-values associated to every single potential change point. Then we eliminate false alarms, {\it i.e. } the change points which have p-value smaller than a given critical level. Next, we apply our method and Penalized Least Squares Criterion procedure to detect change points on simulated data sets and then we compare them. Eventually, we apply the Filtered Derivative with p-Value method to the segmentation of heartbeat time series, and the detection of change points in the average daily volume of financial time series.} \\



\section{INTRODUCTION}
\label{sec1} In a wide variety of applications including health
and medicine, finance, civil engineering, one models time
dependent systems by a sequence of  random variables described by
a finite number of structural parameters. These structural
parameters can change abruptly and it is relevant to detect the
unknown change points. Both on-line and off-line change point detection have
their own relevance, but in this work we are concerned  with
off-line detection.

Statisticians have studied change point detections since the
1950's and there is a huge literature on this subject, see e.g. the textbooks \citet{Basseville:Nikiforov:1993,Brodsky:Darkhovsky:1993, Csorgo:Horvath:1997,Montgomery:1997}, the interesting papers of \citet{Chen:1988,Steinebach:Eastwood:1995}, and let us also refer to \citet{Huskova:Meintanis:2006,Kirch:2008,Gombay:Serban:2009} for update reviews, \citet{Fliess:etal:2010} who propose an algebraic approach for change point detection problem,  and \citet{Birge:Massart:2007} for a good summary of the
model selection approach.

Among the popular methods, we find the Penalized Least Squares
Criterion (PLSC). This algorithm is based on the minimization of
the contrast function when the number of change points is known,
see \citet{Bai:Perron:1998}, \citet{Lavielle:Moulines:2000}.
When the number of changes is unknown, many authors use the penalized
version of the contrast function, see e.g. \citet{Lavielle:Teyssiere:2006} or \citet{Lebarbier:2005}.
From a numerical point of view, the least squares methods are based on dynamic programming
algorithm which needs to compute a matrix. Therefore, the time and
memory complexity 
are of order
$\mathcal{O}(n^2)$ where $n$ is the size of data sets. So,
complexity becomes an important limitation with technological
progress.

Indeed,  recent measurement methods allow us to record and to
stock large data sets.  For example, in Section~\ref{sec5}, we
present change point analysis of heartbeat time series: It is
presently possible to record the duration of  each single
heartbeat during a marathon race or for  healthy people during 24
hours. This leads to data sets of size $n\ge 40,000$ or  $n\ge
100,000$, respectively. Actually, this phenomenon is general: time
dependent data are often recorded at very high frequency (VHF),
which combined with size reduction of memory capacities allows
recording of millions of data.

This technological progress leads us to revisit change point
detection methods in the particular case of large or huge data
sets. This framework constitutes the  main novelty of this work:
we have to develop embeddable algorithms with low time and memory
complexity. Moreover, we can adopt an asymptotic point of view. The non asymptotic case will be treated elsewhere.

For the off-line detection problem, many papers propose procedures based on the dynamic programming algorithm. Among them, as previously said, we find the Penalized Least Squares
Criterion (PLSC) which relies on the minimization of contrast function. It performs with time and space complexities of order 
$\mathcal{O}\left(n^2\right)$, see e.g. \citet{Bai:Perron:1998,Lavielle:Teyssiere:2006} or \citet{Lebarbier:2005}. Then, \citet{Guedon:2007} has developed a method based on the maximization of the log-likelihood. Using a forward-backward dynamic programming algorithm, he obtains a time complexity of order $\mathcal{O}\left((K+1)n^2\right)$ and a memory complexity of order $\mathcal{O}\left((K+1)n\right)$, where $K$ represents the number of change points.\\
These methods are then time consuming when the number of observations $n$ becomes large. For this reason, many other procedures 
have been  proposed for a faster segmentation. In \citet{Chopin:2007} or \citet{Fearnhead:Liu:2007}, the problem of multiple change-points is reformulated as a Bayesian model. The segmentation is computed by using a particle filtering algorithm which requires complexity of order $\mathcal{O}\left(M \times n\right)$ both in time and memory, where $M$ is the number of particles. For instance, \citet{Chopin:2007} obtains 
reasonable estimations of change points with $M=5000$ particles. Recently, a segmentation method based on wavelet decomposition in Haar basis and thresholding has been proposed by \citet{BenYaacov:Eldar:2008}, with time complexity of order  $\mathcal{O}\left(n \log(n)\right)$  and memory complexity of order  $\mathcal{O}\left(n\right)$.
Though, no mathematical proofs has been given. Moreover, no results has been proposed to detect change points in the variance or in the slope and the intercept of linear regression model.

In this paper, we investigate the properties of  a new off-line
detection methods for multiple change points, so-called  Filtered
Derivative with p-value method (FDp-V). Filtered Derivative  has
been introduced by \citet{Benveniste:Basseville:1984}. \citet{Basseville:Nikiforov:1993},
next \citet{Antoch:Huskova:1994} propose an asymptotic study and 
\citet{Bertrand:2000} gives some non asymptotic
results. On the one hand, the advantage of Filtered Derivative method
is its time and memory complexity, both of order $\mathcal{O}(n)$.
On the other hand, the drawback of Filtered Derivative method is
that if it detects the right change points it also gives many
false alarms. To avoid this drawback, we introduce a second step
in order to disentangle right change points and false alarms. In
this second step we calculate the p-value associated to each
potential change point detected in the first step. Stress that the
second step has still time and memory complexity  of order
$\mathcal{O}(n)$.

Our belief is that FDp-V method is quite general for large
datasets. However, in this work, we restrict ourselves to
detection of change points on mean and variance for a sequence of
independent  random variables and change point on slope and
intercept for linear model. The rest of this paper is organized as
follows:   In Section~\ref{sec2}, we describe Filtered Derivative
with p-value Algorithm. Then,  Section~\ref{sec3} is concerned
with theoretical results and FDp-V method for detecting  changes
on  mean and   variance. In Section~\ref{sec4}, we present
theoretical  results of FDp-V method for detecting changes on
slope and intercept for   linear regression model. Finally, in
Section~\ref{sec5}, we give numerical simulations with a
comparison with PLSC algorithm, and we present
some results on real data showing the robustness of our method (as no independance is guaranteed).
All the proofs are postponed to Section~\ref{sec6}.


\section{DESCRIPTION OF THE FILTERED DERIVATIVE WITH p-VALUE
METHOD}  \label{sec2}

In this section, we describe the Filtered Derivative with p-value
method (FDp-V). First, we describe precisely the statistical model we will use throughout our work. Next, we
describe the two steps of FDp-V method: Step~1 is based on
Filtered Derivative and select the potential change points,
whereas Step~2 calculate the p-value associated to each potential
change point, for disentangling right change points and false
alarms.
\subsubsection*{Our model:}
Let $\displaystyle (X_t)_{t=1,\dots,n}$ be a sequence of
independent r.v. with distribution $\mathcal{M}_{\theta(t)}$ ,
where  $\theta\in \R^d$ is a finite dimensional parameter. We
assume that the maps $t\mapsto \theta(t)$ is piecewise constant,
{\em i.e.} there exists a configuration of change points
$0=\tau_0<\tau_1<\dots<\tau_K<\tau_{K+1}=n$ such that $\theta(t)=
\theta_k$ for $\tau_k \le t < \tau_{k+1}$. The integer $K$
corresponds to the number of change times and $(K+1)$ to the
number of segments. In summary, if $j\in[\tau_k,\tau_{k+1}[$, the r.v. $X_j$ are independent and
identically distributed with distribution
$\mathcal{M}_{\theta_k}$.\\
We stress that the number of abrupt changes $K$ is unknown, leading to
a problem of model selection.
There is a huge literature on change point
analysis and model selection see e.g. the monographs
\citet{Basseville:Nikiforov:1993,Brodsky:Darkhovsky:1993} or \citet{Birge:Massart:2007}. As pointed out in the introduction, large and huge datasets lead to revisit change point analysis by taking into account  time and memory complexity of the different methods.

\subsubsection*{Filtered Derivative:}
Filtered Derivative is defined as the difference between
the estimators of the parameter $\theta$ computed on two sliding
windows respectively at the right and at the left of the index
$k$, both of size $A$, that is specified by the  following function: 
\begin{equation}
D(k,A)  = \widehat{\theta}(k,A)-\widehat{\theta}(k-A,A),
\label{def:D}
\end{equation}
where $\displaystyle \hat{\theta}(k,A)$ is an
estimator of $\theta$ on the sliding box $[k+1,\,k+A]$.
Eventually, this method consists on filtering data by computing
the estimators of the parameter $\theta$ before  applying a
discrete derivation. So, this construction explains the name of the
algorithm, so-called
Filtered Derivative method.

The norm of the filtered derivative function, say $\|D\|$,
exhibits "hats" at the vicinity of parameter  change points.
To simplify the presentation, without losing generality, we will assume now that
$\theta$ is a one dimensional parameter, say for example the mean. Thus, the change points can be estimated as
arguments of local maxima of $|D|$, see Figure~\ref{fig1} below.
Moreover, the size of changes in $\theta$ is equal to the height
of positive or negative hats of $D$.
\begin{figure}[htbp]
\begin{center}
\begin{tabular}{c}
\includegraphics[width=12cm,height=6.5cm]{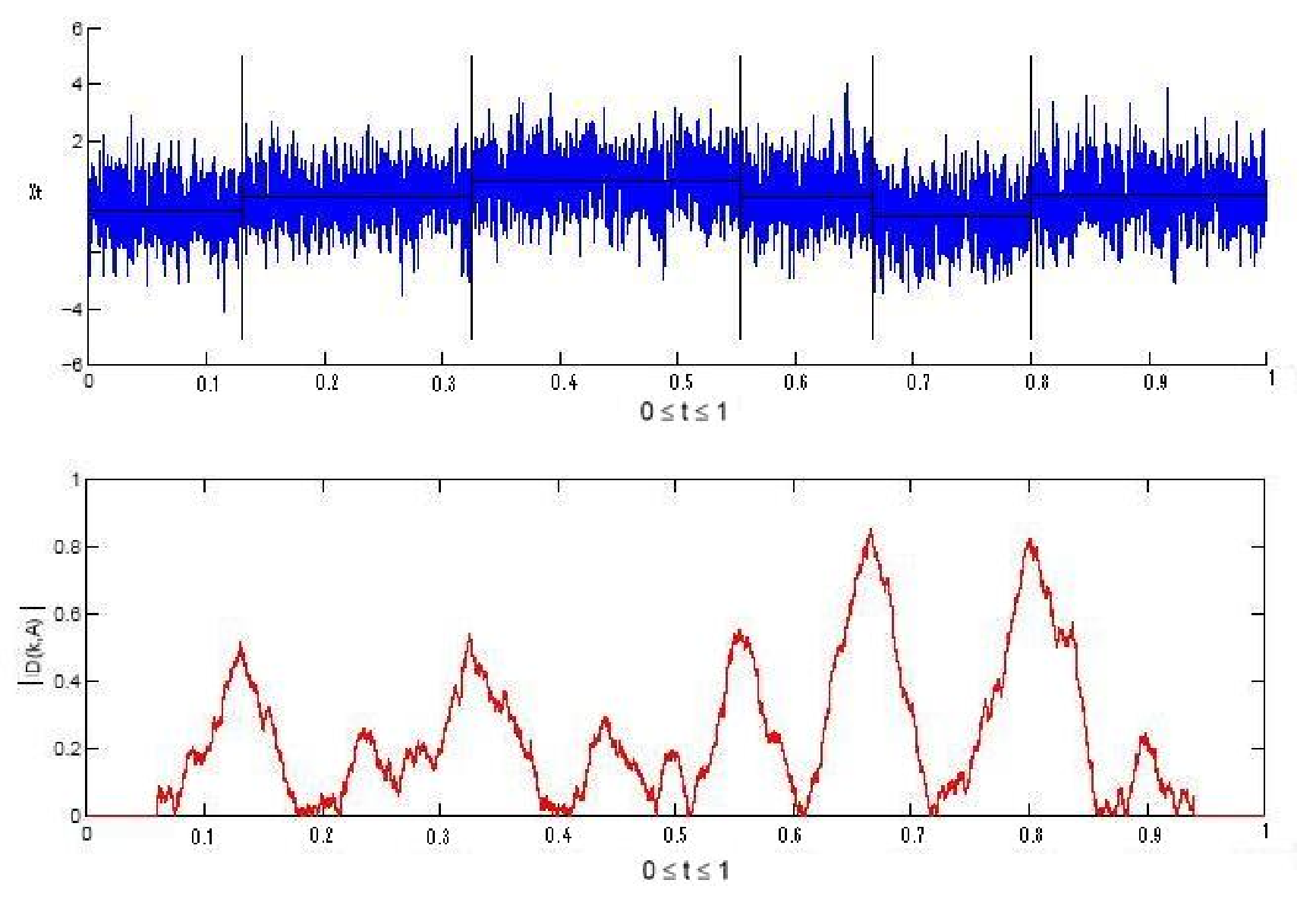}
\end{tabular}
\caption{\emph{Above: Gaussian random variables simulated with
constant variance and presenting changes in the mean. Below:
Filtered derivative function D.}} \label{fig1}
\end{center}
\end{figure}

Nevertheless, we remark through the graph of the function $|D|$ that there are not only the "right hats" (surrounded in blue in Figure~\ref{fig2}) which gives the right change points, but also  false alarms (surrounded in green in Figure~\ref{fig2}). Consequently, we have introduced another idea in order to keep just the right change points. This objective is reached by splitting the detection procedure into two successive steps: In Step~1, we detect  potential change points as local maxima of the filtered derivative function. In Step~2, we test
wether a potential change point is a false alarm or not. Both
steps use estimation of the p-value of existence of change point.

The construction of two different statistical tests and  the computation of p-values is detailed below.

\subsubsection*{Step~1: Detection of the potential change points}
In order to detect the potential change points, we test the null hypothesis $(H_0)$ of no change in the parameter
$\theta$
$$(H_0): \    \theta_1=\theta_2=\cdots=\theta_{n-1}=\theta_n$$
against the alternative hypothesis $(H_1)$ indicating the existence of one or more than one 
changes 
$$(H_1): \    \text{There is an integer } K \in \N^{*} \text{ and }
0=\tau_0<\tau_1<\dots<\tau_K<\tau_{K+1}=n \text{ such that }$$
$$
\theta_1 = \cdots = \theta_{\tau_{1}} \neq \theta_{\tau_{1}+1} =
\cdots = \theta_{\tau_{2}} \cdots \neq \theta_{\tau_{K}+1} =
\cdots = \theta_{\tau_{K+1}}.
$$
where $\theta_i$ is the value of the parameter $\theta$ for $1 \le
i \le n$.

In \citet{
Antoch:Huskova:1994} or \citet{Bertrand:2000},
potential change points are selected as  times $\tau_k$ corresponding to local maxima of 
the absolute value of the filtered derivative $|D(\tau_k,A)|$, when moreover this last quantity 
exceeds a given threshold $\lambda$. However, the efficiency of the 
approach is strongly linked to  the choice of the threshold $\lambda$.
Therefore, in this work, we have a slightly different approach: 
we fix a probability of type I error at level $p_1^*$,
and we determine the corresponding critical value $C_1$ given 
by $$\PP \left( \max_{k \in [A:n-A]} |D(k,A)| >C_1| \ H_0 \text{ is true} \right)=p_1^*.$$ 
Of course, such a probability is usually not available,
so that we only have the asymptotic distribution of the maximum of $|D|$,
which will be the main part of Section~\ref{sec3} and Section~\ref{sec4}.
Then, roughly speaking, we select as potential change points,
local maxima  for which $\ds|D(\tau_k,A)|>C_1$.

By doing so, the p-value  $p_1^*$ appears as an hyper parameter. The second hyper parameter is the window size $A$. As pointed out in
\citet{Antoch:Huskova:1994,Bertrand:2000}, Filtered Derivative
method works under the implicit assumption that the minimal
distance between two successive change points is greater than
$2A$. An adaptive and automatic choice of $A$ is an open and particularly relevant question, but there is no such result in the litterature. Thus, we need some {\it a priori} knowledge on
the minimal length between two successive changes.

More formally, we have the following algorithm:
\\~\\
\textbf{Step~1 of the algorithm}
\begin{enumerate}
\item \textbf{Choice of the hyper parameters}
    \begin{itemize}
    \item Choice of the window size $A$ from information of  the
    practitioners.
    \item Choice of $p_1^*$: \\First we fix the significance level of type 	I error at $p_1^*$. Then, the theoretic expression of type I error, 		given in Section~\ref{sec3} and Section~\ref{sec4}, fixes the value 		of the threshold $C_1$.\\
The list of potential change points, which contains right 			change points as well as false ones, can be too large. So the level 	of significance for Step~1 can be chosen large, \textit{i.e.} $p_1^*=5\% \text{ or } 10\%$.
    \end{itemize}

\item \textbf{Computation of the filtered derivative function}
\begin{itemize}
    \item
The memory complexity results from the recording of the filtered
derivative sequence $\ds \left(D(k,A)\right)_{A \leq k \leq n-A}$. 
Clearly,  we need a second buffer of size $(n-2A)$. Thus the  memory complexity is of order $n$. 
On the other hand, filtered derivative function can be calculated
by recurrence, see $\eqref{rec1}$ or $\eqref{rec2}$ and
$\eqref{def:hat:b}$.
These recurrence formulas induce  time complexity of order $n \times \mathcal{C}$ where $ \mathcal{C}$ is the time complexity of an iteration. For instance, in  $\eqref{rec1}$
we have $ \mathcal{C} \simeq 4\, additions$, in $\eqref{rec2}$ we have  $ \mathcal{C} \simeq 5\, additions + 5\, multiplications + 1\, division$, and 
in $\eqref{def:hat:b}$,  $ \mathcal{C} \simeq 3\, additions + 1\, division$.
\end{itemize}
\item \textbf{Determination of the potential change points}
\begin{itemize}
\item Initialization:~ \\
Set counter of potential change point $k=0$ and
$\widetilde{\tau}_k= \argmax_{k\in [A, n-A]} |D(k,A)|$. \item
While ($\ds|D(\widetilde{\tau}_k,A)|>C_1$) do
\begin{itemize}
\item $k=k+1$
\item $D(k, A) = 0$ for all $k\in( \widetilde{\tau}_{k}-A,\widetilde{\tau}_{k}+A)$.
\item
$\widetilde{\tau}_k= \argmax_{k\in [A, n-A]} |D(k,A)|$
\end{itemize}
We increment the change point counter and we set the values of the function $D$  to zero because the width of the hat is  equal to $2A$.
\item Finally, we sort the vector  $\displaystyle \left( \widetilde{\tau}_1,\ldots,\widetilde{\tau}_{K_{\max}} \right)$
in increasing order. The integer $K_{\max}$ represents the number of potential change points detected at the end of the Step~1. By construction, 
we have $K_{\max} \leq \lfloor n/A\rfloor-2$, where $\lfloor x \rfloor$ corresponds to the integer part of $x$.
\end{itemize}
\end{enumerate}
\textbf{End of the Step~1 of the algorithm}
\\

\begin{figure}[htbp]
\begin{center}
\begin{tabular}{c}
\includegraphics[width=12cm,height=4cm]{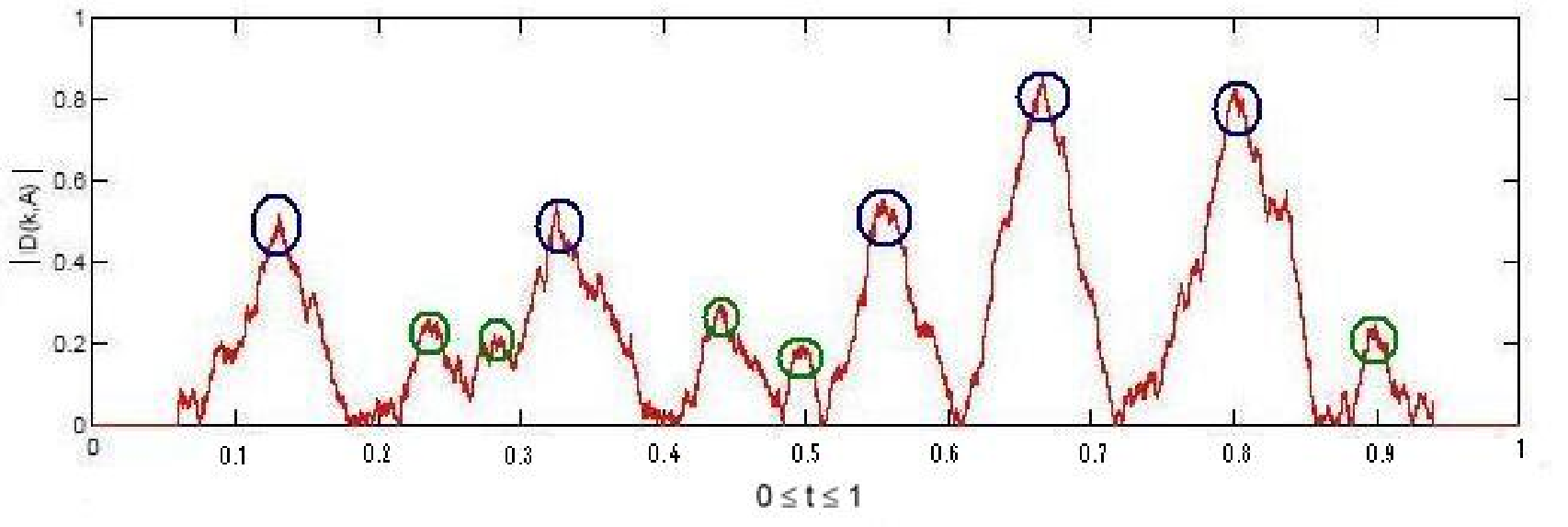}
\end{tabular}
\caption{\emph{Detection of potential change points}} \label{fig2}
\end{center}
\end{figure}

Figure~\ref{fig2} below provides an example: The family of
potential change points contains the right change points
(surrounded in blue in Figure~\ref{fig2}) as well as false alarms
(surrounded in green in Figure~\ref{fig2}). 

\subsubsection*{Step~2: Removing false detection}
The list of potential change points $\displaystyle \left(\widetilde{\tau}_1,\ldots,\widetilde{\tau}_{K_{\max}} \right)$ obtained at step 1 contains right change points but also false detections. 
In the second step a test is carried out to remove the false detection from the list of change points found at step 1.
By doing so, we obtain a subset $\displaystyle \left(\widehat{\tau}_{1},\ldots,\widehat{\tau}_{\widehat{K}} \right)$ of the first list.


%
More precisely, for all potential change point $\widetilde{\tau}_{k}$, we test wether the parameter
is the same on the two successive intervals $(\widetilde{\tau}_{k-1}+1,\widetilde{\tau}_{k})$ and
$(\widetilde{\tau}_{k}+1,\widetilde{\tau}_{k+1})$, or not.
Formally, for all $1 \le k \le K_{\max}$, we apply  the following hypothesis testing
\begin{eqnarray*}
   (H_{0,k}):   \theta_k=  \theta_{k+1}
   &\qquad\text{versus}\qquad&
    (H_{1,k}):
     \theta_k\neq  \theta_{k+1}
     \end{eqnarray*}
where $\theta_k$ is the value of $\theta$ supposed to be constant on the segment $(\widetilde{\tau}_{k-1}+1,\widetilde{\tau}_{k})$. By using this second test, we calculate new p-values $\ds
(\widetilde{p}_{1},\ldots,\widetilde{p}_{K_{\max}})$ associated
respectively to each potential change points $\ds
(\widetilde{\tau}_1,\ldots,\widetilde{\tau}_{K_{\max}})$. Then, we
only keep the change points which have a p-value smaller than
a critical level denoted $p_2^{*}$. Step~2 must be much
more selective with a significance level $p_2^{*}=10^{-3} \text{ or } 10^{-4}$. Consequently, Step~2 permits us to remove most of false detections,  
and so to deduce an estimator of the piecewise constant map $t\mapsto \theta(t)$, see Figure~\ref{fig3} below.\\

\begin{Remark}
We note that $\widetilde{\tau}_{k}$, for $ 1 \leq k \leq K_{\max}$, is a random variable.  So 
that at first sight, the average value of the parameter $\theta$ on the interval $(\widetilde{\tau}_{k-1}+1,\widetilde{\tau}_{k})$ may also be a random variable. However, \citet{Bertrand:2000} has showed that, with probability converging to 1,  the right change point is in $(\widetilde{\tau}_{k}+\epsilon_k,\widetilde{\tau}_{k}-\epsilon_k)$, where $\epsilon_k$ is a given bounded real number. For the detection of change in the mean in a series  of independent Gaussian random variables, \citet[Theorem 3.1, p.255]{Bertrand:2000} has proved that $\epsilon_k=5(\sigma/\delta_k)^2+1$ where $\delta_k$ represents the size of change in the mean at $\tau_k$. This result can be easily generalized to a series of independent non Gaussian random variables with finite second moment. Therefore, on a set $\Omega_n$ with $\PP(\Omega_n )\rightarrow 1$ as $n \rightarrow \infty$, the average value of $\theta$ on the segment $(\widetilde{\tau}_{k}+\epsilon_k,\widetilde{\tau}_{k}-\epsilon_k)$  becomes not random and actually is equal to $\theta_k$.
\end{Remark}

The time and memory complexity of this second step of the
algorithm is still $\mathcal{O}(n)$ because we only need to
compute and to stock $(K_{\max}+1)$ estimators of parameter
$\theta$ which are successively compared two by two in order to
eliminate false alarms.

\begin{figure}[htbp]
\begin{center}
\begin{tabular}{c}
\includegraphics[width=12cm,height=6cm]{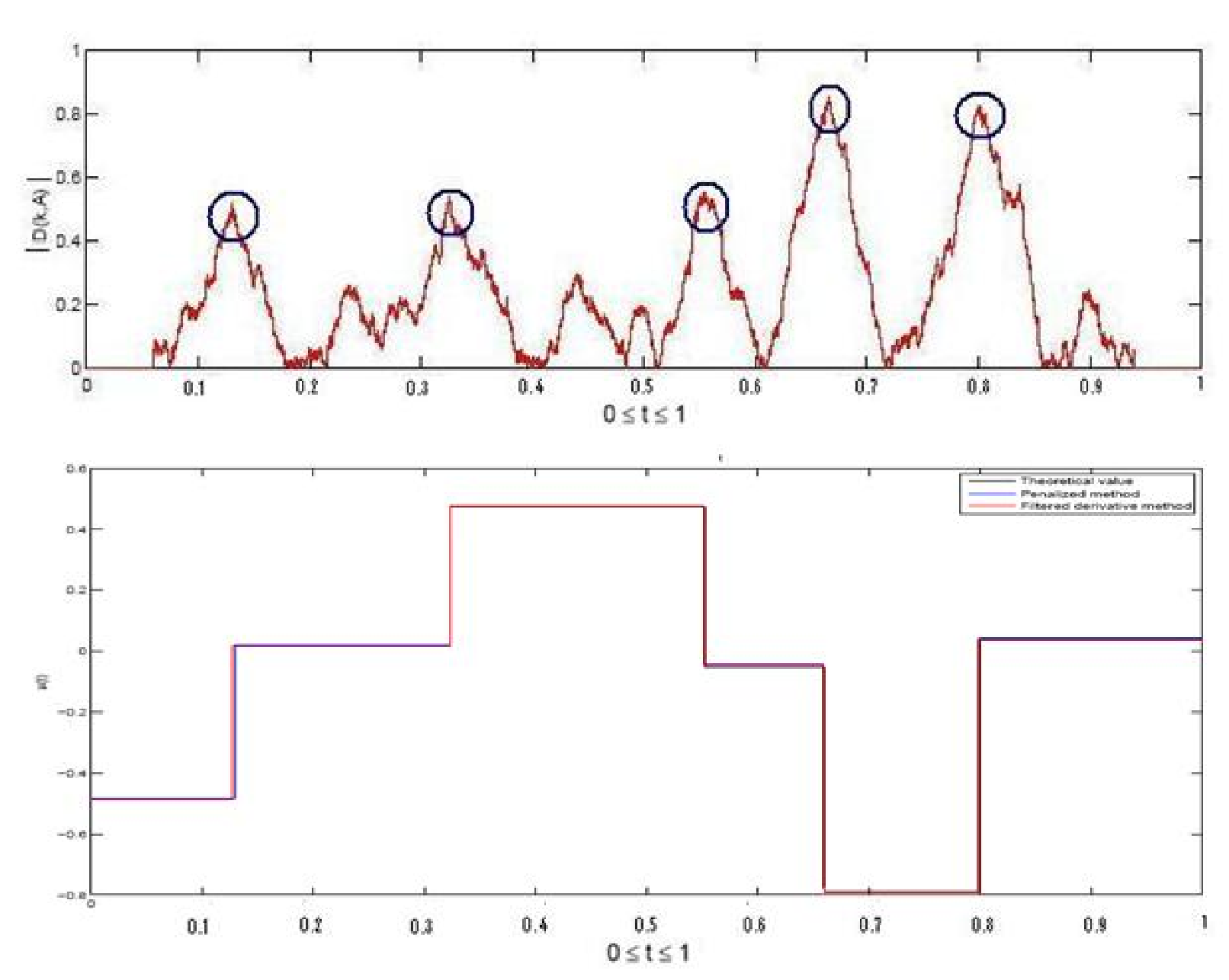}
\end{tabular}
\caption{\emph{Above: Detection of right change points. Below:
Theoretical value of the piecewise-constant map $t\mapsto \mu_t$
(black), and its estimators given by PLSC method (blue) and
FDp-V method (red).}} \label{fig3}
\end{center}
\end{figure}
To sum up, FDp-V method is a two step procedure, in which the  filtered derivative method is applied first and then a second test is carried out to remove the false detection from the list of change points found in Step 1. 
The first step has both time and  memory  complexity of order  $\mathcal{O}(n)$. At 
the second step, the number of selected potential change points   is  smaller. As  a consequence, both time and  memory  complexity of Step 2 are still of order  $\mathcal{O}(n)$.
 
Let us finish the presentation of the FDpV algorithm by the following remark: in this paper we will restrict ourselves to the detection of change points of a one dimensional parameter, such as mean, variance, slope and intercept of linear regression. However, our algorithm still works in any finite dimension parameter space. Indeed, for each parameter, we can compute the corresponding filtered derivative sequence and its asymptotic distribution and then compare it to the  corresponding type I errors~$p^*_1$ and $p^*_2$.


\section{THEORETICAL RESULTS}
\label{sec3}
In this section, we present theoretical results of the FDp-V method for detecting changes on the mean and on the variance. The two following subsections give the asymptotic distribution of the filtered derivative used for the detection of potential change points. Then, Subsection~\ref{subsec:pvalue} provides formulae for computing the p-values during Step~2 in order to eliminate false alarms. 

\subsection{Change in the mean}
Let $\displaystyle (X_i)_{i=1,\dots,n}$ be a sequence of
independent r.v. with mean $\mu_i$ and a known
variance~$\sigma^2$. We assume that the map $i\mapsto \mu_i$ is
piecewise constant, {\em i.e.} there exists a configuration
$0=\tau_0<\tau_1<\dots<\tau_K<\tau_{K+1}=n$ such that $\E(X_i)=
\mu_k$ for $\tau_k \le i < \tau_{k+1}$. The integer $K$
corresponds to the number of changes. However, in any real life
situation, the number of abrupt changes $K$ is unknown, leading to
a problem of model selection, see e.g. \citet{Birge:Massart:2007}.
\text{}\\
Filtered Derivative method applied to the mean is based on the
difference between the empirical mean computed on two sliding
windows respectively at the right and at the left of the index
$k$, both of size $A$, see
\citet{Antoch:Huskova:1994,Basseville:Nikiforov:1993}. 
This difference corresponds to a sequence $(D_1(k,A))_{A \leq k \leq
n-A}$ defined by \ba\label{def:D1}
D_1(k,A)&=&\hat{\mu}(k,A)-\hat{\mu}(k-A,A) \ea
where
$ \hat{\mu}(k,A)=\frac{1}{A} \sum_{j=k+1}^{k+A}X_j$
is the empirical mean of $X$ on the (sliding) box $[k+1,\,k+A]$.
These quantities can easily be calculated by recurrence with
complexity $\mathcal{O}(n)$. It suffices to remark  that
\begin{equation}
 A\times D_1(k+1,A)=A\times D_1(k,A)+X_{k+A+1}-2X_k+X_{k-A+1}. \label{rec1}
\end{equation}
First we give in Theorem~\ref{th1} the asymptotic
behaviour of the maximum of $|D_1|$ under null hypothesis of no
change in the mean and with size of the sliding windows tending to
infinity at a certain rate. In the sequel, we denote by $A_n$ the
size of the sliding windows and we will always suppose that \ba\label{An}
\lim\limits_{\substack{n \to +\infty}} \frac{A_n}{n}=0 \text{ and
} \lim\limits_{\substack{n \to +\infty}} \frac{ (\log
n)^2}{A_n}=0. \ea

\begin{theo}[Change point in the mean with known variance]
\label{th1} Let $\displaystyle (X_i)_{i=1,\dots,n}$ be a sequence
of independent identically distributed random variables with mean
$\mu$, variance $\sigma^2$ and assume that one of the following
assumptions is satisfied
\begin{description}
    \item [$\mathbf{(\mathcal{A}_1)}$] $\ X_1 \sim \mathcal{N}(\mu,\sigma^2).$
    \item [$\mathbf{(\mathcal{A}_2)}$] $\exists\, t > 0 \text{ such as } \ \E[\exp(tX_1)] < +\infty
    \text{, and }
    \displaystyle
    \lim\limits_{\substack{n \to +\infty}}
    \frac{(\log n)^3}{A_n}=0.$
    \item [$\mathbf{(\mathcal{A}_3)}$] $\exists\, p > 2 \text{ such as } \E[|X_1|^p] < +\infty
    \text{, and }
    \displaystyle \lim\limits_{\substack{n \to +\infty}}
    \frac{n^{2/p} \log n}{A_n}=0.$
\end{description}
Let $D_1$ be defined by (\ref{def:D1}). Then under the null hypothesis
\begin{equation}
\lim\limits_{\substack{n \to +\infty}} \PP \left( \max_{k \in
[A_n:n-A_n]} |D_1(k,A_n)| \le \frac{\sigma}{\sqrt{A_n}} c_n(x)
\right) = \exp(-2e^{-x}), \label{lim:mu}
\end{equation}
\begin{equation}
 c_n(x)=c \left ( \frac{n}{A_n} -1,x \right )\label{cnx}
\end{equation}
\begin{equation}
c(y,x)=\frac{1}{\sqrt{2\log y}} \left ( x +  2 \log y + \frac 1 2
\log \log y -\frac 1 2 \log \pi \right ). \label{cyx}
\end{equation}
\begin{flushright}
$\square$
\end{flushright}
\end{theo}
A similar version of Theorem~\ref{th1} was first proved by \citet{Chen:1988}, via the strong invariance principle, in order to detect changes in the mean. Later, \citet{Steinebach:Eastwood:1995} have studied the same results for the detection of changes in the intensity of the renewal counting processes, and have improved them by giving another rate for the size of the sliding window under assumptions $\mathbf{(\mathcal{A}_2)}$ and $\mathbf{(\mathcal{A}_3)}$. The reader is also referred to the book of \citet{Revesz:1990} for a good summary about the increments on partial sums.\\
\text{}\\
Furthermore, for a known change point $\tau_k$, the law of the random variable $D_1(\tau_k,A_n)$ defined by $\eqref{def:D1}$ is known. Therefore, the probability to detect a known change point in the mean is described in the following remark.  

\begin{Remark} [Probability to detect known change point in the mean] 
Let $\tau_k$ be a known change point in the mean of size $\delta_k:=|\mu_{\tau_k}-\mu_{\tau_k+1}|$. Then, under assumption $\mathbf{(\mathcal{A}_1)}$, the probability to detect it is given by
\ba \label{proba:detect:1a} \PP \left( |D_1(\tau_k,A_n)| \geq C_1\right) = 1-\phi \left(\frac{(C_1-\delta_k)\sqrt{A_n}}{\sqrt{2}\sigma} \right)+\phi \left(\frac{(-C_1-\delta_k)\sqrt{A_n}}{\sqrt{2}\sigma} \right), \ea and under assumptions $\mathbf{(\mathcal{A}_2)}$ and $\mathbf{(\mathcal{A}_3)}$, we have the following asymptotic probability
\ba \label{proba:detect:1b} \lim\limits_{\substack{n \to +\infty}}  \PP \left( \sqrt{A_n} |D_1(\tau_k,A_n)-\delta_k| \geq C_1\right) = 2 \phi \left(\frac{-C_1}{\sqrt{2}\sigma} \right) \ea
where $C_1$ is the threshold fixed in Step~1 and $\phi$ is the standard normal distribution. The proof can be deduced from $\ds D_1(\tau_k,A_n)=\delta_k+\frac{\sqrt{2}\sigma}{\sqrt{A_n}} Z $ where $Z$ is the standard normal r.v.
\end{Remark} 

In applications, the variance $\sigma^2$ is unknown. For this
reason we may replace it by its empirical estimator,
$\widehat{\sigma}_{n}^2$. But, in order to keep the same result as in Theorem~\ref{th1}, the estimator $\widehat{\sigma}_{n}^2$ has
to verify a certain condition given by the following theorem.

\begin{theo}[Change point in the mean with unknown variance]
\label{th2} We apply to the same notations and the same
assumptions as in Theorem~\ref{th1}. Moreover, we assume that
$\widehat{\sigma}_{n}$ is an estimator of $\sigma$ satisfying
\be\label{lim:hyp1} \lim\limits_{\substack{n \to +\infty}}
|\sigma-\widehat{\sigma}_{n}| \log n \stackrel{P}{=} 0, \ee
where the sign $\stackrel{P}{=}$ means convergence in probability.
Then, under the null hypothesis,
\begin{equation}
\lim\limits_{\substack{n \to +\infty}} \PP \left( \max_{k \in
[A_n:n-A_n]} |D_1(k,A_n)| \le
\frac{\widehat{\sigma}_n}{\sqrt{A_n}} c_n(x) \right) =
\exp(-2e^{-x}). \label{lim:mu2}
\end{equation}
\begin{flushright}
$\square$
\end{flushright}
\end{theo}
Remark that condition (\ref{lim:hyp1}) is not really restrictive, indeed as soon as $\widehat\sigma^2_n$ satisfies a CLT, the condition is verified. For example, with the usual empirical variance estimator, a fourth order moment (of $X$) is sufficient.
\subsection{Change in the variance}
Now, we consider the case where we have a set of observations
$\displaystyle (X_i)_{i=1,\dots,n}$ and we wish to know whether
their variance has changed at an unknown time. If $\mu$ is known,
then the problem is very simple. Testing $(H_0)$ against $(H_1)$
means that we are looking for a change in the mean of the sequence
$\displaystyle \left( (X_i-\mu)^2 \right )_{i=1,\dots,n}$.\\
Filtered Derivative method applied to the variance is based on the
difference between the empirical variance computed on two sliding
windows respectively at the right and at the left of the index
$k$, both of size $A_n$ which satisfy condition $\eqref{An}$. This
difference is in fact a sequence of random variables denoted by
$(D_2(k,A_n))_{A_n \leq k \leq n-A_n}$ and defined as follows \ba
\label{def:D2}
D_2(k,A_n)&=&\hat{\sigma}^2(k,A_n)-\hat{\sigma}^2(k-A_n,A_n) \ea
where
\ban
\hat{\sigma}^2(k,A_n)&=&\frac{1}{A_n}
\sum_{j=k+1}^{k+A_n}(X_j-\mu)^2
\ean
denotes the empirical variance of $X$
on the box $[k+1,k+A_n]$.
By using Theorem \ref{th1}, we can deduce 
 the asymptotic distribution of the maximum of
$|D_2|$ under null
hypothesis $(H_0)$. This gives straightforwardly the following corollary:
\begin{coro}[Change point in the variance with known mean]
\label{cor1} Let $\displaystyle (X_i)_{i=1,\dots,n}$ be a sequence
of independent identically distributed random variables with mean
$\mu$ and assume that one of the following assumptions is
satisfied
\begin{description}
    \item [$\mathbf{(\mathcal{A}_4)}$] $\exists\, t > 0 \text{ such as } \ \E[\exp(t(X_1-\mu)^2)] < +\infty
    \text{, and }
    \displaystyle
    \lim\limits_{\substack{n \to +\infty}}
    \frac{(\log n)^3}{A_n}=0.$
    \item [$\mathbf{(\mathcal{A}_5)}$] $\exists\, p > 2 \text{ such as } \ \E \left[|X_1-\mu|^{2p} \right ] < +\infty
    \text{, and }
    \displaystyle \lim\limits_{\substack{n \to +\infty}}
    \frac{n^{2/p} \log n}{A_n}=0.$
\end{description}
Let $D_2$ be defined by $\eqref{def:D2}$ and $\nu^2=\var \left
[(X_1-\mu)^2 \right ]$. Then, under the null hypothesis,
\begin{equation}
\lim\limits_{\substack{n \to +\infty}} \PP \left( \max_{k \in
[A_n:n-A_n]} |D_2(k,A_n)| \le \frac{\nu}{\sqrt{A_n}} c_n(x)
\right) = \exp(-2e^{-x}),
\end{equation}
where $c_n(.)$ is defined by $\eqref{cnx}$.\hfill
$\square$
\end{coro}
In practical situations we rarely know the value of the
constant mean. However, $\hat{\mu}(k,A_n)$ is a consistent estimator on the box
$[k+1,\,k+A_n]$ for the mean  $\mu$ under both null hypothesis $(H_0)$ and alternative hypothesis $(H_1)$.
Thus, in the definition of the sequence $D_2$, $\mu$ can be  replaced by its
estimator $\hat{\mu}(k,A_n)$ on the box $[k+1,\,k+A_n]$. To be precise, 
let \ba\label{def:D2hat}
\widehat{D}_2(k,A_n)&=&\widetilde{\sigma}(k,A_n)-\widetilde{\sigma}(k-A_n,A_n), \ea
where 
\ban \widetilde{\sigma}(k,A_n)&=&\frac{1}{A_n}
\sum_{j=k+1}^{k+A_n}(X_j-\hat{\mu}(k,A_n))^2
\ean
denotes the empirical
variance of $X$ on the sliding box $[k+1,\,k+A_n]$ with unknown
mean. So, in order to obtain the same asymptotic distribution as
the one obtained in Corollary~\ref{cor1}, we must add extra
conditions on the estimators $\hat{\mu}(k,A_n)$. These new
conditions are given in the next corollary.
\begin{coro}[Change point in the variance with unknown mean]
\label{cor2} With the notations and assumptions of
Corollary~\ref{cor1}, we suppose moreover that
\begin{equation}
\lim\limits_{\substack{n \to +\infty}} \max_{0 \leq k \leq
n-A_n}|\mu-\hat{\mu}(k,A_n)| \left ( A_n \log n \right
)^{\frac{1}{4}} \stackrel{a.s}{=} 0 \label{condMu}
\end{equation}
where the sign $\stackrel{a.s}{=}$ means almost surely
convergence. Then under the null hypothesis
\begin{equation}
\lim\limits_{\substack{n \to +\infty}} \PP \left( \max_{k \in
[A_n:n-A_n]} |\widehat{D}_2(k,A_n)| \le \frac{\nu}{\sqrt{A_n}}
c_n(x)  \right) = \exp(-2e^{-x}). \label{lim:var2}
\end{equation}
\begin{flushright}
$\square$
\end{flushright}
\end{coro}
Let us remark that (\ref{condMu}) is not a very stringent condition.

\subsection{Step~2: calculus of p-values for changes on the mean and changes on the variance}
\label{subsec:pvalue}
In this subsection, we recall p-value formula associated with the second test in order to remove false alarms. Let us stress that the only
novelty of this subsection is the idea to divide the detection of abrupt
change into two steps, see section \ref{sec3}. Since the r.v. $X_i$ are independent,
the calculus of p-value relies on well known results that can be found in any statistical textbook.

First, consider the Gaussian case and let us introduce some
notations: For $1 \leq k \leq K_{\max}$, let
$(X_{\widetilde{\tau}_{k-1}+1},\ldots,X_{\widetilde{\tau}_{k}})$
and
$(X_{\widetilde{\tau}_{k}+1},\ldots,X_{\widetilde{\tau}_{k+1}})$
two successive samples of i.i.d. Gaussian random variables such that \ban X_{\widetilde{\tau}_{k-1}+1} \sim
\mathcal{N}(\mu_k,\sigma_k^2)& \quad\text{ and }\quad&
X_{\widetilde{\tau}_{k}+1} \sim
\mathcal{N}(\mu_{k+1},\sigma_{k+1}^2). \ean We can use Fisher's F-statistic to determine the p-value of the existence of a change on the variance  at time $\widetilde{\tau}_{k}$ under the null assumption (H0):
$\ds \sigma_k^2=\sigma_{k+1}^2$. Then, we can use Student statistic to determine the p-value of a change on the mean  at time $\widetilde{\tau}_{k}$, that is under the null assumption
 (H0): $\ds \mu_k=\mu_{k+1}$.

Second, consider the general case. Since $A_n\to\infty$ and by construction $\ds \widetilde{\tau}_{k+1}-\widetilde{\tau}_k\ge A_n$, we can apply CLT as soon as the r.v. $X_i$ satisfy Lindeberg condition.  Thus, the empirical  mean and the empirical variance converge to the corresponding ones in the Gaussian case.  Eventually, if one of the assumptions $(\mathcal{A}_i)$ for $i=1$, $2$ or $3$,  then we can still apply Fisher and Student statistics to compute the p-value.


 \section{LINEAR REGRESSION}
\label{sec4}
In this Section, we consider linear regression model:
\ba \label{ordYi}
Y_i &=& a X_i + b + e_i, \text{ for }
1 \le i \le n,
\ea
where the terms $\left( e_i
\right)_{i=1,\dots,n}$ are independent and identically distributed
Gaussian random errors with zero-mean and variance $\sigma^2$ and, to simplify,
$\left( X_i
\right)_{i=1,\dots,n}$ are  equidistant time points given by
\ba\label{absXi}
 X_i=i\Delta&\quad \mathrm{with}\quad& \Delta > 0.
\ea
Our aim is to detect change points on the parameters $(a, b)$ of the linear model.
As the Filtered Derivative is a local method, and to simplify notations, we will restrict ourselves here to one change point.

The two following subsections give the result for detection
of potential change points on the slope and on the intercept.
Subsection \ref{subsect:pvalue:linear:regression}
provides formulas for calculating the p-values during Step~2.

\subsection{Change in the slope}
In this subsection, we are concerned with detection of change points in
the slope $a$ when the intercept remains constant. 

Filtered Derivative method applied to the slope is based on the
differences between   estimated values of the slope
$a$ computed on two sliding windows  at the right and
at the left of the index $k$, both of size $A_n$.
These differences, for $k\in [A_n, n-A_n]$, form a sequence of
random variables,  given by
 \ba \label{def:D3}
D_3(k,A_n)&=&\hat{a}(k,A_n)-\hat{a}(k-A_n,A_n) \ea
where
\begin{equation}
\hat{a}(k,A_n)=\left[ A \sum_{j=k+1}^{k+A}X_j
Y_j-\sum_{j=k+1}^{k+A}X_j \sum_{j=k+1}^{k+A}Y_j\right] \,\left[  A
\sum_{j=k+1}^{k+A}X_j^2 - \left (\sum_{j=k+1}^{k+A}X_j \right
)^2\right]^{-1}, \label{rec2}
\end{equation}
is the estimator of the slope
$a$ on the (sliding) box $[k+1,\,k+A_n]$.  Let us stress that these quantities can  be
calculated by recurrence with complexity $\mathcal{O}(n)$.\\
Our first result gives the asymptotic distribution of the maximum
of $|D_3|$ under the null hypothesis of no change on the linear regression.

\begin{theo}[Change point in the slope]
\label{th3}
Let $\left( X_i
\right)_{i=1,\dots,n}$ and $\left( Y_i \right)_{i=1,\dots,n}$ be given by
$\eqref{absXi}$ and $\eqref{ordYi}$  where $e_i$ is a family of i.i.d. mean zero Gaussian r.v. with variance $\sigma^2$.
Let $D_3$ be defined by $\eqref{def:D3}$ and assume that $A_n$
satisfies condition $\eqref{An}$. Then under the null hypothesis
\begin{equation}
\lim\limits_{\substack{n \to +\infty}} \PP \left( \max_{k \in
[A_n:n-A_n]} |D_3(k,A_n)| \le \frac{2\sqrt{6}\sigma}{\Delta
\sqrt{A_n(A_n^2-1)}} d_n(x) \right) =
\exp(-2e^{-x}), \label{lim:a}
\end{equation}
with
$d_n(x)=c \left ( A_n,x \right )
$ and
 $c(.,.)$ is defined by $\eqref{cyx}$. \hfill $\square$
\end{theo}

Moreover, for a given change point $\tau_k$, the law of the random variable $D_3(\tau_k,A_n)$ defined by $\eqref{def:D3}$ is known. Therefore, the probability to detect a known change point in the slope is specified in the following remark.

\begin{Remark} [Probability to detect known change point in the slope] Let $\tau_k$ be a known change point in the slope of size $\delta_{a,k}:=|a_{\tau_k}-a_{\tau_k+1}|$. Then the probability to detect it is given by
\ba \label{proba:detect:3} \PP \left( |D_3(\tau_k,A_n)| \geq C_1\right) = 1-\phi \left(\frac{(C_1-\delta_{a,k})\Delta \sqrt{A_n(A_n^2-1)}}{2\sqrt{6}\sigma} \right)\nonumber\\\qquad\qquad\,\,\,+\phi \left(\frac{(-C_1-\delta_{a,k})\Delta \sqrt{A_n(A_n^2-1)}}{2\sqrt{6}\sigma} \right) \ea where $C_1$ is the threshold fixed in Step~1 and $\phi$ is the standard normal distribution. The proof is directly deduced from $ \ds D_3(\tau_k,A_n)=\delta_{a,k}+\frac{2\sqrt{6}\sigma}{\Delta \sqrt{A_n(A_n^2-1)}} Z $ where $Z$ is a standard normal random variable.
\end{Remark} 

\subsection{Change in the intercept}
In this subsection, we focus on the detection of change points in
the intercept with known
slope $a$. To do this, we calculate the differences between  estimators of
the intercept $b$ computed on two sliding windows respectively at
the right and at the left of the index $k$, both of size $A_n$.
For
$A_n \leq k \leq n-A_n$,
these differences form
a sequence of random variables given by
\ba \label{def:D4}
D_4(k,A_n)&=&\hat{b}(k,A_n)-\hat{b}(k-A_n,A_n)
\ea
where
\ba\label{def:hat:b}
 \hat{b}(k,A_n)&=&\frac{1}{A} \sum_{j=k+1}^{k+A}Y_j - a
\times \frac{1}{A} \sum_{j=k+1}^{k+A}X_j,
\ea
is the estimator of
the intercept $b$ on the (sliding) box $[k+1,\,k+A_n]$.
By applying Theorem~\ref{th1}, we get
the asymptotic distribution of the maximum
of $|D_4|$ under the null hypothesis of no change on the linear regression.

\begin{coro}[Change point in the intercept] \label{cor3}
With the notations and
assumptions of Theorem \ref{th3}, we have under the null hypothesis
\begin{equation}
\lim\limits_{\substack{n \to +\infty}} \PP \left( \max_{k \in
[A_n:n-A_n]} |D_4(k,A_n)| \le \frac{\sigma}{\sqrt{A_n}} c_n(x)\right) = \exp(-2e^{-x}), \label{lim:b}
\end{equation}
where $c_n(.)$ is defined by $\eqref{cnx}$.\hfill $\square$
\end{coro}
\text{}\\
In addition, for a known change point $\tau_k$, the distribution of the random variable $D_4(\tau_k,A_n)$ defined by $\eqref{def:D4}$ is known. Therefore, the probability to detect a known change point in the intercept is given in the following remark.

\begin{Remark} [Probability to detect known change point in the intercept] Let $\tau_k$ be a known change point in the intercept of size $\delta_{b,k}:=|b_{\tau_k}-b_{\tau_k+1}|$. Then the probability to detect it is given by
\ba \label{proba:detect:4} \PP \left( |D_4(\tau_k,A_n)| \geq C_1\right) = 1-\phi \left(\frac{(C_1-\delta_{b,k})\sqrt{A_n}}{\sqrt{2}\sigma} \right)+\phi \left(\frac{(-C_1-\delta_{b,k})\sqrt{A_n}}{\sqrt{2}\sigma} \right) \ea where $C_1$ is the threshold fixed in Step~1 and $\phi$ is the standard normal distribution. The proof may be easily deduced from $\ds D_4(\tau_k,A_n)=\delta_{b,k}+\frac{\sqrt{2}\sigma}{\sqrt{A_n}} Z $ where $Z$ is a standard normal random variable.
\end{Remark} 

In real applications the slope is often unknown. In this case, we replace it  in $\eqref{def:D4}$
by its empirical estimator
$\widehat{a}_n$. This leads to the definition
 \ba \label{def:D4hat}
\widehat{D}_4(k,A_n)&=&\widetilde{b}(k,A_n)-\widetilde{b}(k-A_n,A_n)
 \ea
where $$ \widetilde{b}(k,A_n)=\frac{1}{A} \sum_{j=k+1}^{k+A}Y_j -
\widehat{a}_n \times \frac{1}{A} \sum_{j=k+1}^{k+A}X_j,$$ is the
estimator of the intercept $b$ on the (sliding) box $[k+1,\,k+A_n]$ with unknown slope.\\
Naturally, we must assume that the estimator of the slope satisfy
a certain convergence condition which is given in the following
corollary.
\begin{coro}[Change point in the intercept with unknown slope]
\label{cor4} Under the same notations and the same
assumptions than in  Corollary~\ref{cor3}. Moreover, we assume that
the estimator $\hat{a}_n$ of $a$ satisfies the following condition
\begin{equation}
\lim\limits_{\substack{n \to +\infty}}|a-\hat{a}_n|
A_n^{\frac{3}{2}} \Delta_n \sqrt{\log n}  \stackrel{a.s}{=} 0
\label{condA}
\end{equation}
where the sign $\stackrel{a.s}{=}$ means almost surely
convergence. Then under the null hypothesis
\begin{equation}
\lim\limits_{\substack{n \to +\infty}} \PP \left( \max_{k \in
[A_n:n-A_n]} |\widehat{D}_4(k,A_n)| \le \frac{\sigma}{\sqrt{A_n}}
c_n(x)  \right) = \exp(-2e^{-x}).
\label{lim:b2}
\end{equation}
\end{coro}
\begin{flushright}
$\square$
\end{flushright}

\subsection{Step~2: calculus of p-values}
\label{subsect:pvalue:linear:regression} Let us give here
p-value formulae associated to Step~2 in linear regression model $\eqref{absXi}$ and $\eqref{ordYi}$. More precisely, we are concerned with detection of right change points on the slope or intercept. We recall that Step~2 of FDp-V method has been introduced in order to eliminate false alarms. Indeed, at Step~1, a time $X_{\widetilde\tau_k}$ has been selected as a potential change
point.
At Step~2,  we test whether the slopes and intercepts of two data sets at left and right of the potential change point $X_{\widetilde\tau_k}$ are
significantly different or not, and we measure this by the corresponding p-value.

Before going further, let us introduce some notations.   For $1 \leq k \leq K_{\max}$, let
$\mathcal{R}_k= \left
((X_{\widetilde{\tau}_{k-1}+1},Y_{\widetilde{\tau}_{k-1}+1}),\ldots,
(X_{\widetilde{\tau}_{k}},Y_{\widetilde{\tau}_{k}})\right )$ and
$\mathcal{R}_{k+1}= \left
((X_{\widetilde{\tau}_{k}+1},Y_{\widetilde{\tau}_{k}+1}),\ldots,
(X_{\widetilde{\tau}_{k+1}},Y_{\widetilde{\tau}_{k+1}})\right )$
two successive samples
of observations such that the relationship between variables $X$
and $Y$ is given by $\eqref{absXi}$ and $\eqref{ordYi}$, or more explicitly by
\ban
Y_{j}&=&a_k X_{j} + b_k + e_{j} \qquad\qquad\text{ for }\, \widetilde{\tau}_{k-1}+1  \leq j \leq
\widetilde{\tau}_{k}.
\ean

By using definition of the error terms,we can  easily see that slope estimator $\widehat{a}_k$ and intercept estimator $\widehat{b}_k$
 have Gaussian distribution given by
\begin{eqnarray*}
\widehat{a}_k  & \sim & \mathcal{N} \left (a_k,  \sigma^2_{a_k}
\right )
\text{ with } \frac{\sigma_{a_k}^2}{\sigma^2} = \left( \sum_{j=\widetilde{\tau}_{k-1}+1}^{\widetilde{\tau}_{k}} (X_j-\overline{X}_k)^2\right)^{-1}, \\
\widehat{b}_k  & \sim & \mathcal{N} \left (b_k, \sigma^2_{b_k}
\right ) \text{ with } \frac{\sigma_{b_k}^2}{\sigma^2} = \frac{1}{n_k} + \overline{X}_k\left(
\sum_{j=\widetilde{\tau}_{k-1}+1}^{\widetilde{\tau}_{k}}
(X_j-\overline{X}_k)^2\right)^{-1},
\end{eqnarray*}
respectively, where $\overline{X}_k$ is the empirical mean of the sequence $\ds
\left(X_{\widetilde{\tau}_{k-1}+1},\ldots,X_{\widetilde{\tau}_{k}}\right)$.
In the sequel, we denote by $\widehat{\sigma}_{a_k}^2$ and
$\widehat{\sigma}_{b_k}^2$ the empirical variance of respectively
the random variables $\widehat{a}_k$ and $\widehat{b}_k$.
\subsubsection*{Comparing slope}
We want to test if the samples $\mathcal{R}_k$ and
$\mathcal{R}_{k+1}$
present a change in slope or not.\\
\ban
 \left(H_{0,k}^a \right): \   a_k= a_{k+1}
&\quad\text{ against }\quad& \left(H_{1,k}^a \right): \
a_k \ne a_{k+1}. \ean Then, the
p-value $\widetilde{p}_{k,a}$ associated to the potential change
point $\widetilde{\tau}_{k}$ in order to eliminate false alarms
for changes in slope is given by
$$
\widetilde{p}_{k,a}=1-\phi_a \left (
\frac{|\widehat{a}_k-\widehat{a}_{k+1}|}{\sqrt{\frac{\widehat{\sigma}^2_{a_k}}{n_k}+\frac{\widehat{\sigma}^2_{a_{k+1}}}{n_{k+1}}}}
\right)
$$
\text{}\\
where $\phi_a$ is a Student T-distribution with $\upsilon_a$
degrees of freedom
$$\upsilon_a = N \left(\widehat{\sigma}^2_{a_k},\widehat{\sigma}^2_{a_{k+1}},n_k,n_{k+1} \right)=
\left\lfloor \frac{ \left
(\frac{\widehat{\sigma}^2_{a_k}}{n_k}+\frac{\widehat{\sigma}^2_{a_{k+1}}}{n_{k+1}}
\right)^2} { \left (\frac{\widehat{\sigma}^2_{a_k}}{ n_k
\sqrt{n_k-1}} \right )^2 + \left (
\frac{\widehat{\sigma}^2_{a_{k+1}}}{n_{k+1} \sqrt{n_{k+1}-1}}
\right )^2 } \right\rfloor$$ and $\lfloor x\rfloor$ is the integer part of $x$.

\subsubsection*{Comparing intercept}
We desire to study if the intercept before and after a detected change points $\tilde\tau_k$ are equal or not, i.e. if it is a real change point, resulting in testing 
\ban (H_{0,k}^b) \
b_k=b_{k+1} &\quad\text{ against }\quad&
(H_{1,k}^b) \ b_k \ne b_{k+1}.
\ean 
Then, the p-value
$\widetilde{p}_{k,b}$ associated to the potential change point
$\widetilde{\tau}_{k}$ in order to remove false alarms from the list of 
changes in intercept found in step 1, is given by
$$
\widetilde{p}_{k,b}=1-\phi_b \left (
\frac{|\widehat{b}_k-\widehat{b}_{k+1}|}{\sqrt{\frac{\widehat{\sigma}^2_{b_k}}{n_k}+\frac{\widehat{\sigma}^2_{b_{k+1}}}{n_{k+1}}}}
\right)
$$
where $\phi_b$ is a Student T-distribution with $\ds \upsilon_b= N
\left(\widehat{\sigma}^2_{b_k},\widehat{\sigma}^2_{b_{k+1}},n_k,n_{k+1}
\right)$ degrees of freedom.


\section{NUMERICAL RESULTS}
\label{sec5}
In this section, we apply FDpV and PLSC methods to detect abrupt changes in the mean of simulated Gaussian r.v and  we compare the results obtained by both methods. Next, we use the FDpV for the detection of change points in the slope of  linear regression model. Finally, we run the FDpV algorithm to the segmentation of heartbeat time series, and  average daily volume drawn from the financial market.

\subsection{A toy model: off-line detection of abrupt changes in
the mean of independent Gaussian random variables with known
variance}

In the following subsection, we consider the elementary problem, namely: The off-line detection of multiple change points in
the mean of  simulated independent Gaussian r.v with known variance, and we numerically compare the
efficiency of the different estimators given by the FDpV and the PLSC procedures.

\subsubsection*{Numerical simulation}
To begin, for $n=5000$ we have simulated a sequence of Gaussian r.v $\left ( X_1,\ldots,X_n \right)$ with known variance $\sigma^2=1$ and mean $\mu_i=g(i/n)$ where
$g$ is a piecewise-constant function with five  change points, \textit{i.e.} we have chosen
a configuration $1=\tau_0<\tau_1<\dots<\tau_K<\tau_{K+1}=n$ with $K=5$ and
such as $\mu_i= g(k)$ for $\tau_k \le i < \tau_{k+1}$ .  The size of the change point at $\tau_k$ specified by  $\delta_k:=|\mu_{\tau_k}-\mu_{\tau_k+1}|$ are to be chosen in the intervalle $[0.5,1.25]$. \\
Then, we have computed the function $k \rightarrow |D_1(A,k)|$ with $A = 300$, see Figure~\ref{fig1}.\\
Both methods, namely Filtered Derivative with p-values, $p_1^{*}=0.05$ and $p_2^{*}=10^{-4}$, and Penalized Least Squares Criterion provide right results, see Figure~\ref{fig3}. This example is plainly confirmed by Monte-Carlo simulations.

\subsubsection*{Monte-Carlo simulation}
\begin{figure}[htbp]
\begin{center}
\begin{tabular}{c}
\includegraphics[width=12cm,height=3cm]{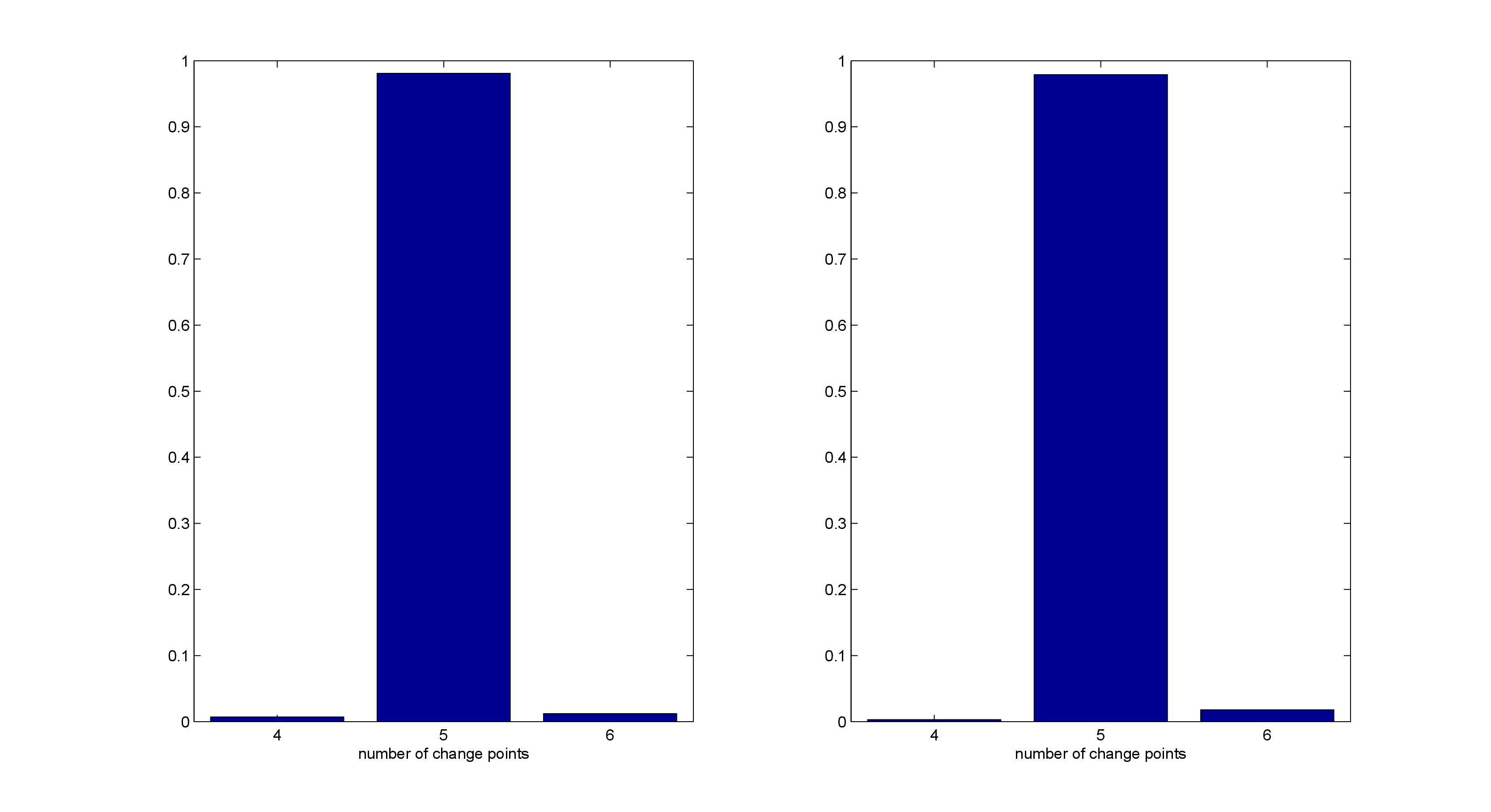}
\end{tabular}
\caption{\emph{Distribution of the estimated number of change
points $\widehat{K}$ for $M=1000$ realizations. Left:~Using PLSC
method. Right: Using Filtered Derivative method.}}\label{fig7}
\end{center}
\end{figure}

In this paragraph, we have made $M=1000$ simulations of
independent copies of sequences of Gaussian r.v.
$X_{1}^{(k)},\ldots,X_{n}^{(k)}$ with variance $\sigma^2=1$ and
mean $\mu(i)=g(i/n)$, for $k=1,\dots,M$. On each sample, we apply
the FDp-V algorithm and the PLSC method. We find the right number of
changes in $98.1\%$ of all cases  for the first method and in
$97.9\%$ for the second one, see Figure~\ref{fig7}.

Then, we compute the mean square errors. There are two kinds of mean square
errors:
\begin{itemize}
  \item
  Mean Integrate Square Error:
$\E \left(\frac{1}{n} \sum_{i=1}^n \left|\widehat{g}(i/n) -
g(i/n)\right|^2 \right) = \E \| \widehat{g} - g
\|_{L^2([0,1])}^{2}$.The estimated function is obtained in two
steps: first we estimate the configuration of change points
$(\hat{\tau}_k)_{k=1,\dots, \hat{K}}$, then we estimate the value
of $\widehat{g}$ between two successive change points as the
empirical mean. \item Square Error on Change Points:
$\displaystyle \E  \left( \sum_{k=1}^K \left|\widehat{\tau}_k -
\tau_k\right|^{2}\right) $, in the case where we have found the
right number of  change points.
\end{itemize}
Table 1 gives the result of Monte Carlo simulation mean errors, and also the comparison between the mean time complexity and the mean
memory complexity. We have written the two programs in Matlab and
have runned it with computer system which has the following
characteristics: 1.8GHz processor and $512$MB memory. \\
\begin{table}[htbp]
\begin{center}
\begin{tabular}{|c|c|c|}
  \hline
   \  & Square Error on Change Points   &  Mean Integrated Squared Error\\
  \hline
  FDp-V method & $1.1840 \times 10^{-4}$ & $0.0107$ \\
  \hline
  PLSC method & $1.2947 \times 10^{-4}$  & $0.0114$ \\
  \hline
  \hline
   \  & Memory allocation (in Megabytes) &  CPU time (in second)\\
  \hline
  FDp-V method & $0.04$ MB & $0.005$ s \\
  \hline
  PLSC method & $200$ MB  & $240$ s \\
  \hline
\end{tabular}
\caption{Errors and complexities given by FDp-V method and PLSC
method.}\label{tab1}
\end{center}
\end{table}

\subsubsection*{Numerical  conclusion}
On the one hand, both methods have the same accuracy in terms of
percentage of the right number of changes, and in terms of Square Error or Mean Integrate Square Error. On the other hand, the
Filtered Derivative with p-Value is less expensive in terms of
time complexity and memory complexity, see Table~\ref{tab1}.
Indeed, Penalized Least Squares Criterion algorithm needs 200 Megabytes 
of computer memory, while Filtered derivative method only  needs $0.008\%$. This plainly confirms the difference of time and memory
complexity, {\em i.e.} $\mathcal{O}(n^2)$ {\em versus} $\mathcal{O}(n)$.

\subsection{Off-line detection of changes in
the slope of simple linear regression}

\begin{figure}[H]
\begin{center}
\includegraphics[width=8.5cm,height=3cm]{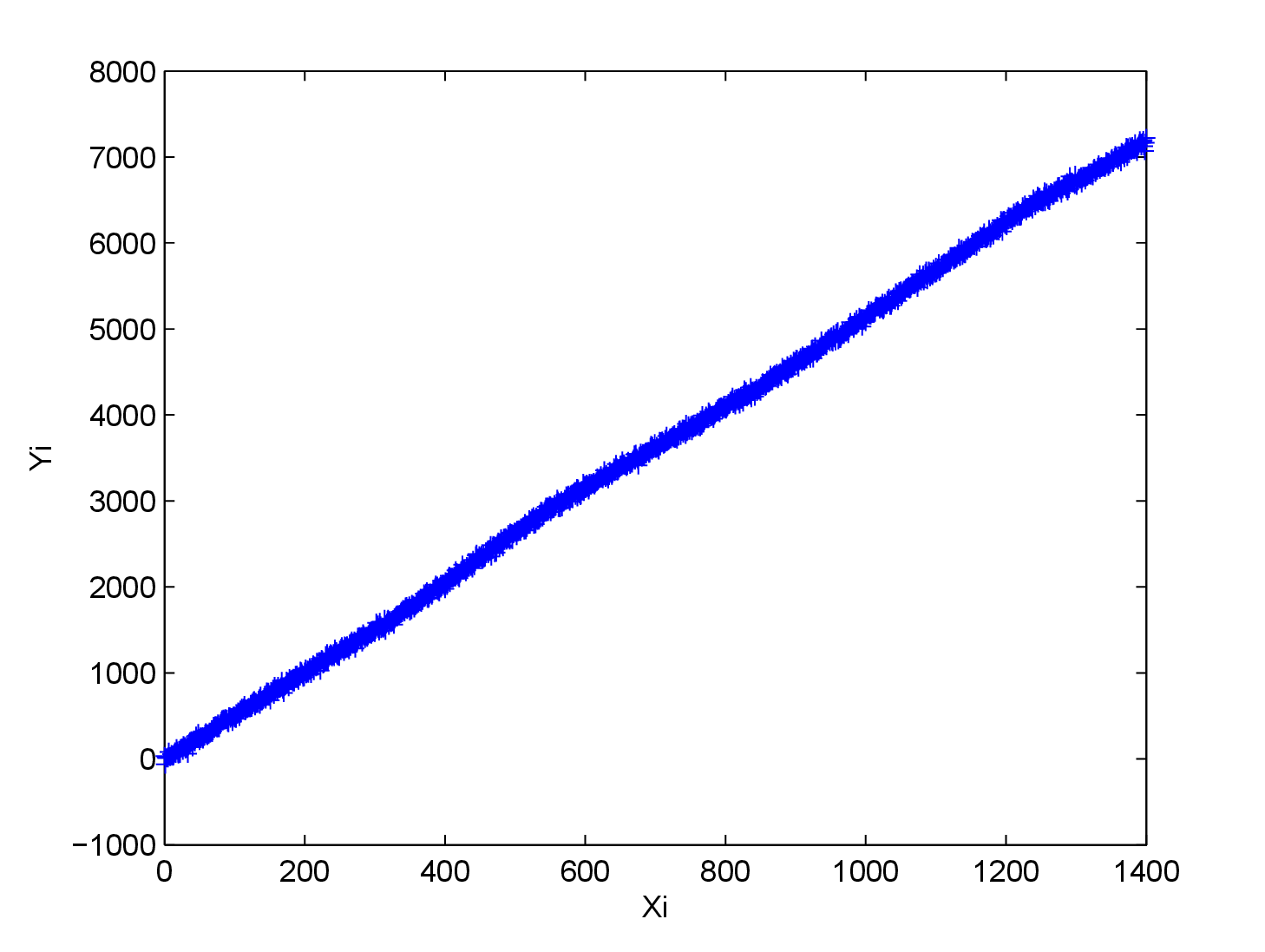}
\caption{\emph{Case 1: Scatter plot of the simulated data
$(X_i,Y_i)$ for $1 \leq i \leq n$.}}\label{fig8}
\end{center}
\end{figure}

\begin{figure}[H]
\begin{center}
\includegraphics[width=8.5cm,height=3cm]{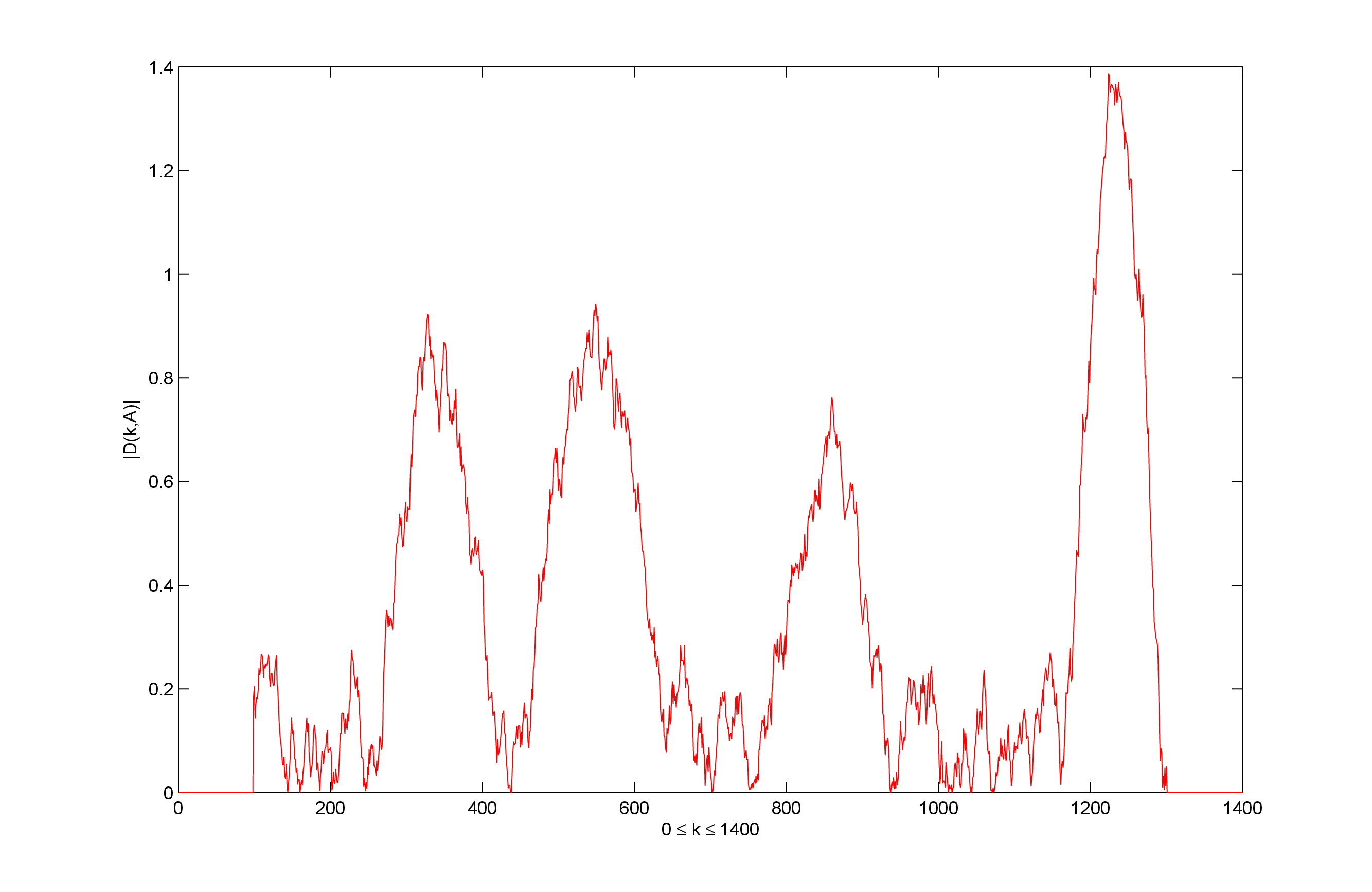}
\caption{\emph{Case 1: The hat function |$D_3$|}.}\label{fig10}
\end{center}
\end{figure}

\begin{figure}[H]
\begin{center}
\includegraphics[width=8.5cm,height=3cm]{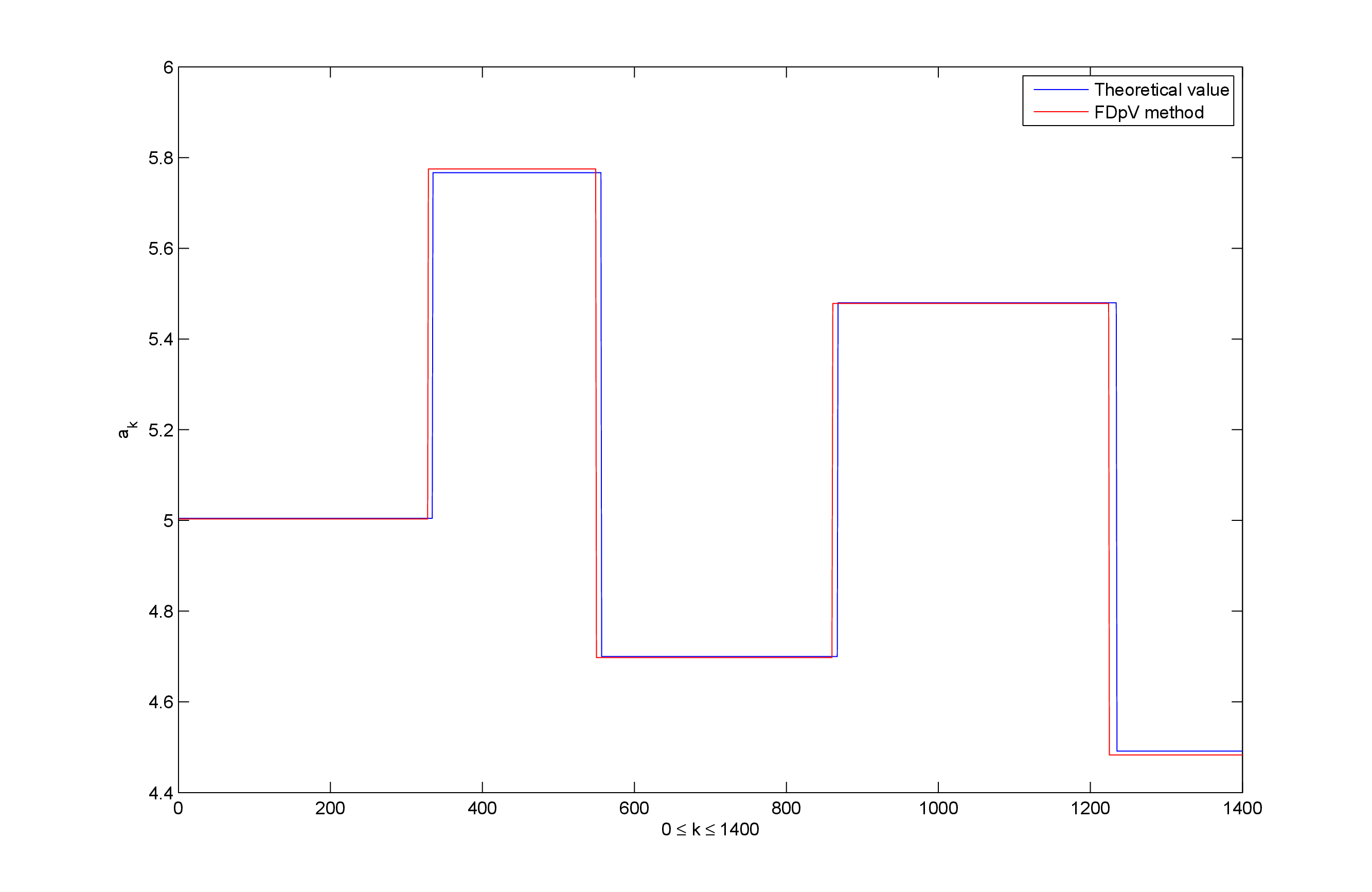}
\caption{\emph{Case 1: Theoretical value of the piecewise-constant
function $h$ (blue), and its estimator given by FDp-V method
(red)}.}\label{fig12}
\end{center}
\end{figure}

\begin{figure}[H]
\begin{center}
\includegraphics[width=8.5cm,height=3cm]{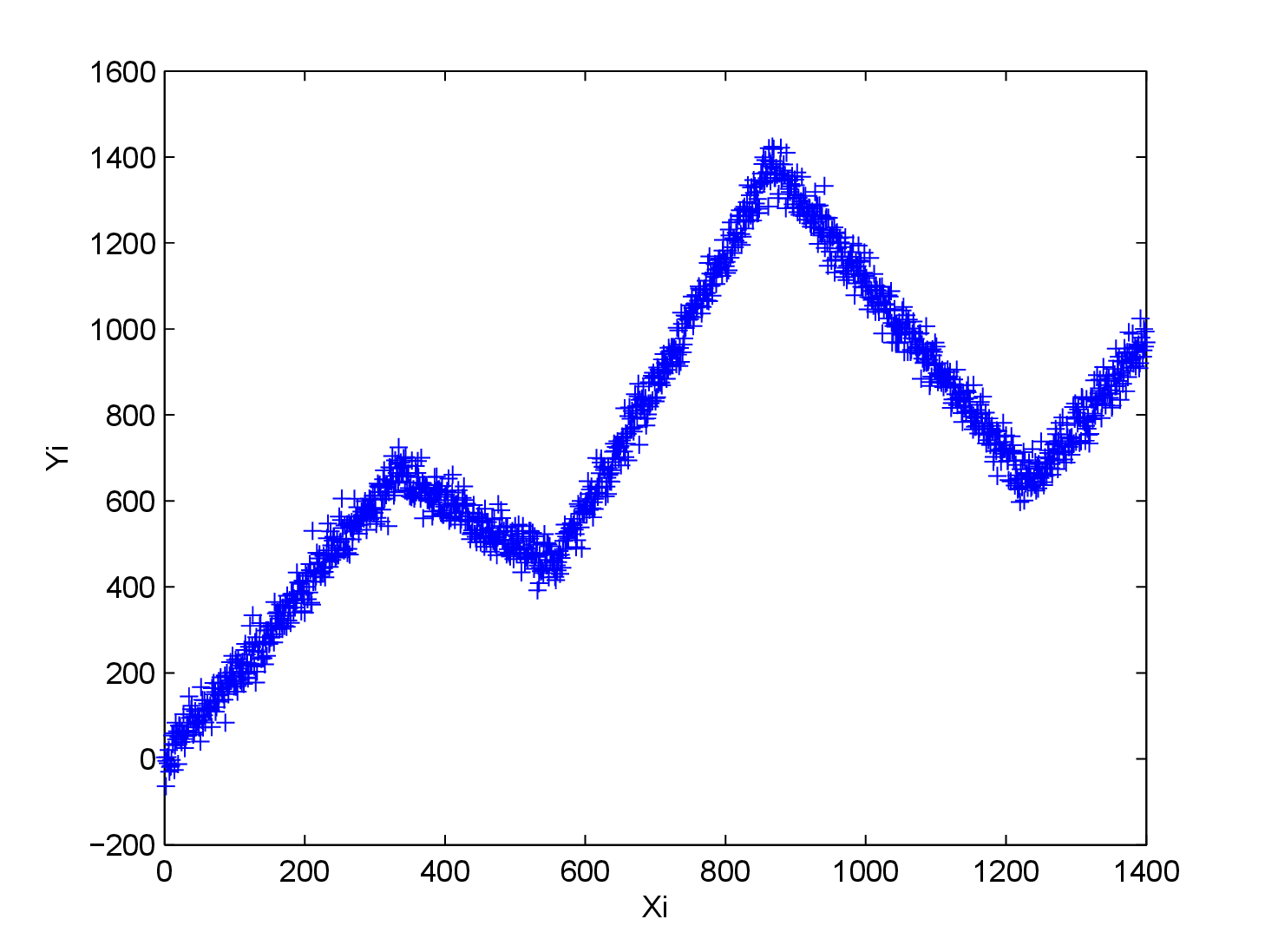}
\caption{\emph{Case 2: Scatter plot of the simulated data
$(X_i,Y_i)$ for $1 \leq i \leq n$.}}\label{fig9}
\end{center}
\end{figure}

\begin{figure}[H]
\begin{center}
\includegraphics[width=8.5cm,height=3cm]{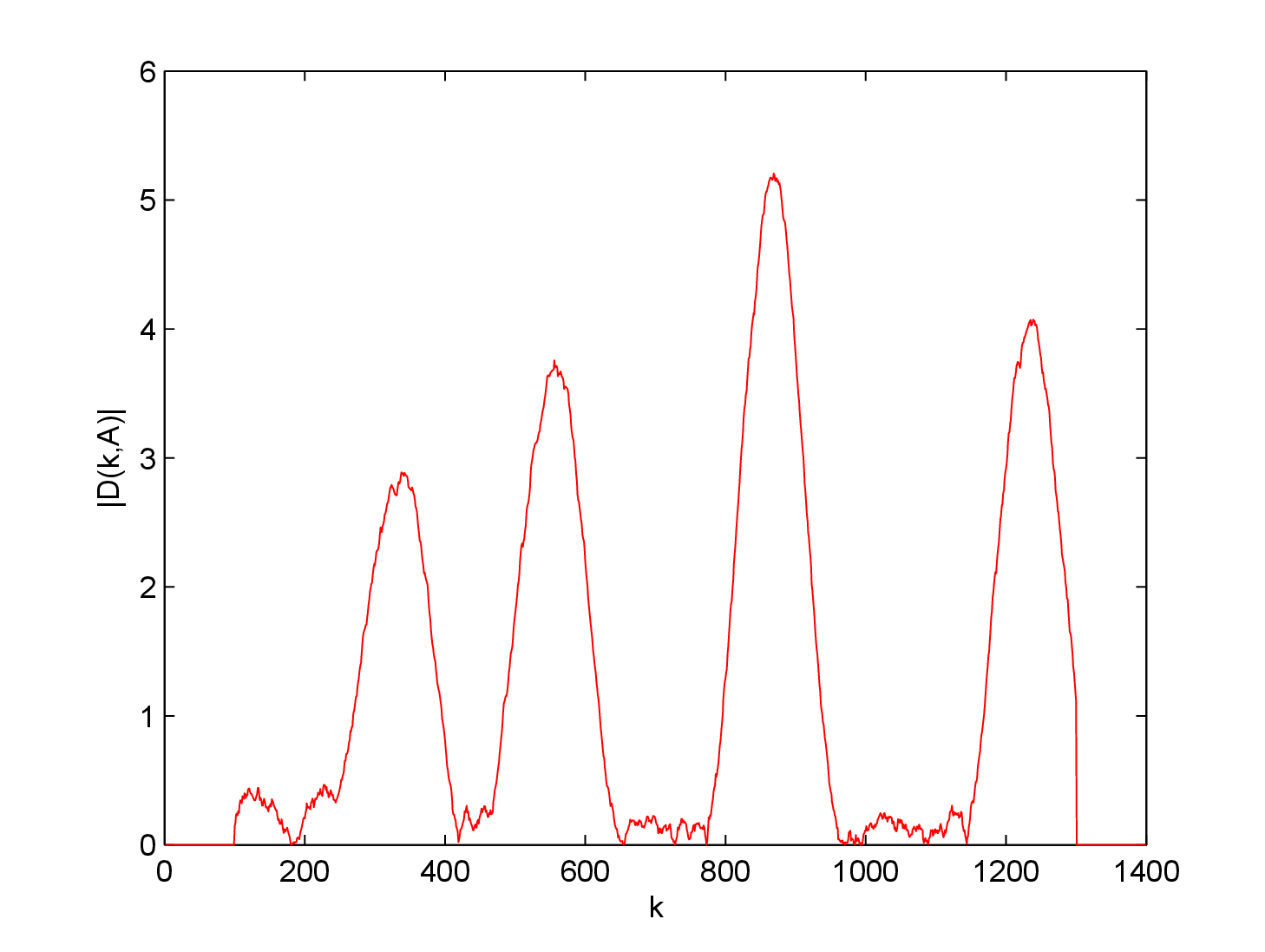}
\caption{\emph{Case 2: The hat function |$D_3$|}.}\label{fig11}
\end{center}
\end{figure}

\begin{figure}[H]
\begin{center}
\includegraphics[width=8.5cm,height=3cm]{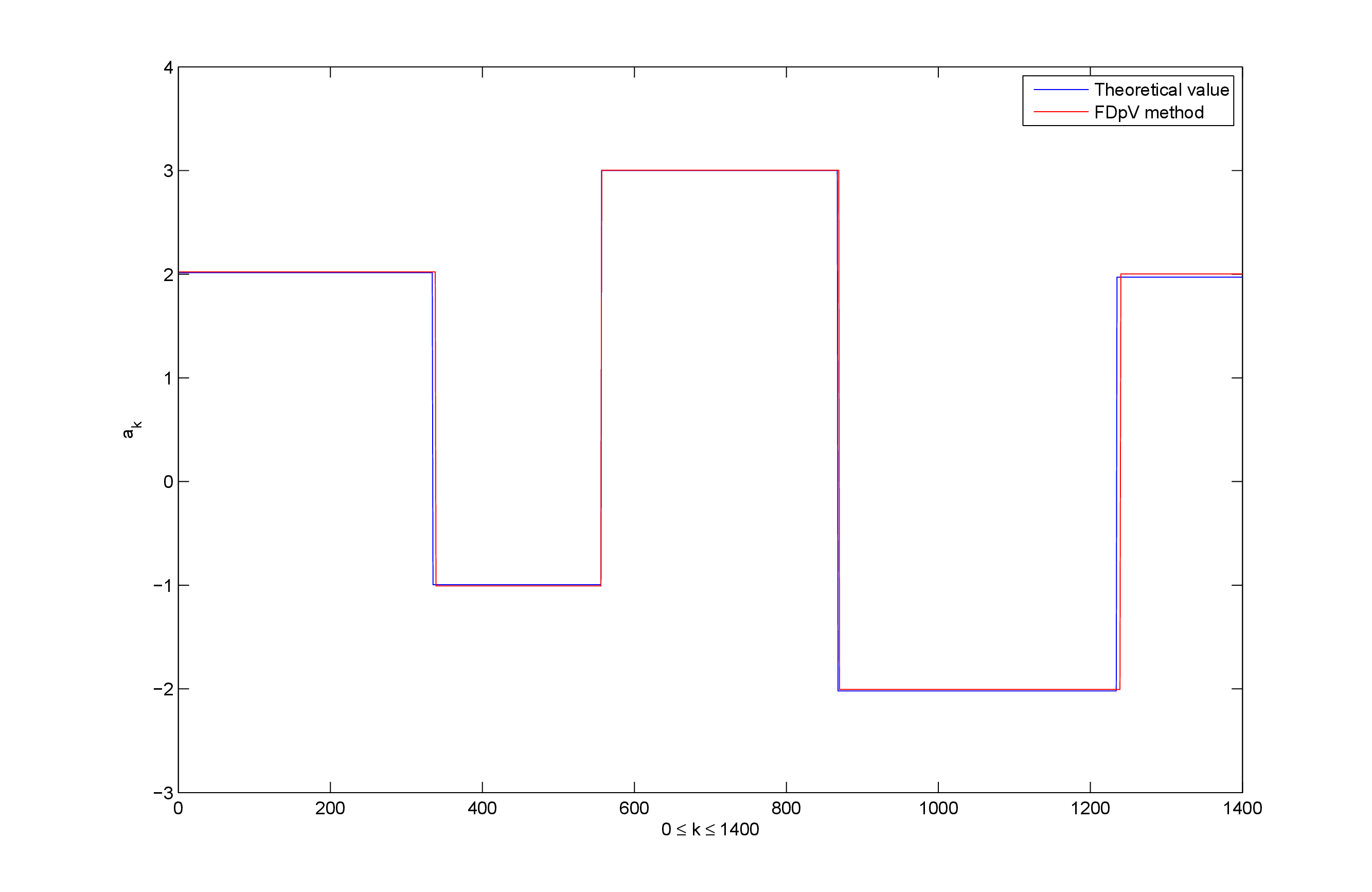}
\caption{\emph{Case 2: Theoretical value of the piecewise-constant
function $h$ (blue), and its estimator given by  FDp-V  method
(red)}.}\label{fig13}
\end{center}
\end{figure}

In this subsection, we consider the problem of 
multiple change points detection in the slope of linear model corrupted by
an additive Gaussian noise. \\
At first, for $n=1400$ we have
simulated the sequences $X=\left ( X_1,\ldots,X_n
\right)$ and $Y=\left ( Y_1,\ldots,Y_n \right)$ defined by
$\eqref{absXi}$ and $\eqref{ordYi}$ with $\Delta=1$, $\sigma=30$
and $a_i=h(i)$ where $h$ is a piecewise-constant function
with four change points, \textit{i.e.} we have chosen
a configuration $1=\tau_0<\tau_1<\dots<\tau_K<\tau_{K+1}=n$ with $K=4$ and
such as $a_i= h(k)$ for $\tau_k \le i < \tau_{k+1}$ . On the one hand, we have considered large change points in the slope such as $\nu_k \in [3,5]$ where
$\nu_k:=|a_{\tau_k}-a_{\tau_k+1}|$ represents the size of the change point at $\tau_k$. On the other hand, we have chosen smaller change points such as $\nu_k \in [0.75,1]$. Next, we have plot $X$ versus
$Y$, see the scatter plots \ref{fig8} and \ref{fig9}. Then, to detect the change points in the slope of these simulated
data, we have computed the function $k \mapsto |D_3(k,A)|$
with $A=100$, see Figures~\ref{fig10}~and~\ref{fig11}.
\text{}\\
Finally, by applying FDp-V procedure with p-values $p_1^{*}=0.05$
and $p_2^{*}=10^{-10}$, we obtain a right localization of the
change points and so a right estimation of the piecewise-constant
function $h$, see Figures~\ref{fig12}~and~\ref{fig13}.

\begin{multicols}{2}

\end{multicols}

\subsection{Application to real data}
In this subsection, we apply our algorithm to detect change
points in the mean of two real samples: the first one is concerned with health and wellbeing, and the second one with finance. Our main purpose here is to show that our estimation procedure is sufficiently robust to consider non-independent time series, that we will study theoretically in a subsequent paper. We moreover provide an analysis of the obtained results.

 \subsection*{An Application to Change Point Detection of Heartbeat Time Series}
In this paragraph, 
 we give examples of application of FDp-V
method to heartbeat time series analysis.   
Electrocardiogram
(ECG)   has been processed from a long time since the implementation of
 monitoring by Holter in  the fifties.
We consider here the RR interval, which provides an accurate  measure of the length of each single heartbeat and corresponds to the instantaneous speed of the heart engine, see \citet{Task:Force:1996}.
From the beginning of 21st century, the size reduction of the measurement
devices allows to record heartbeat time series  for
healthy people in ecological situations throughout  long periods of time:  
Marathon runners,
individuals daily (24 hours) records,  etc.
We then obtain large data sets of more than
40.000 observations, resp. 100,000 observations, 
 that address change detection of heart rate. 
\begin{figure}[htbp]
\begin{center}
\begin{tabular}{c}
\includegraphics[width=16cm,height=6cm]{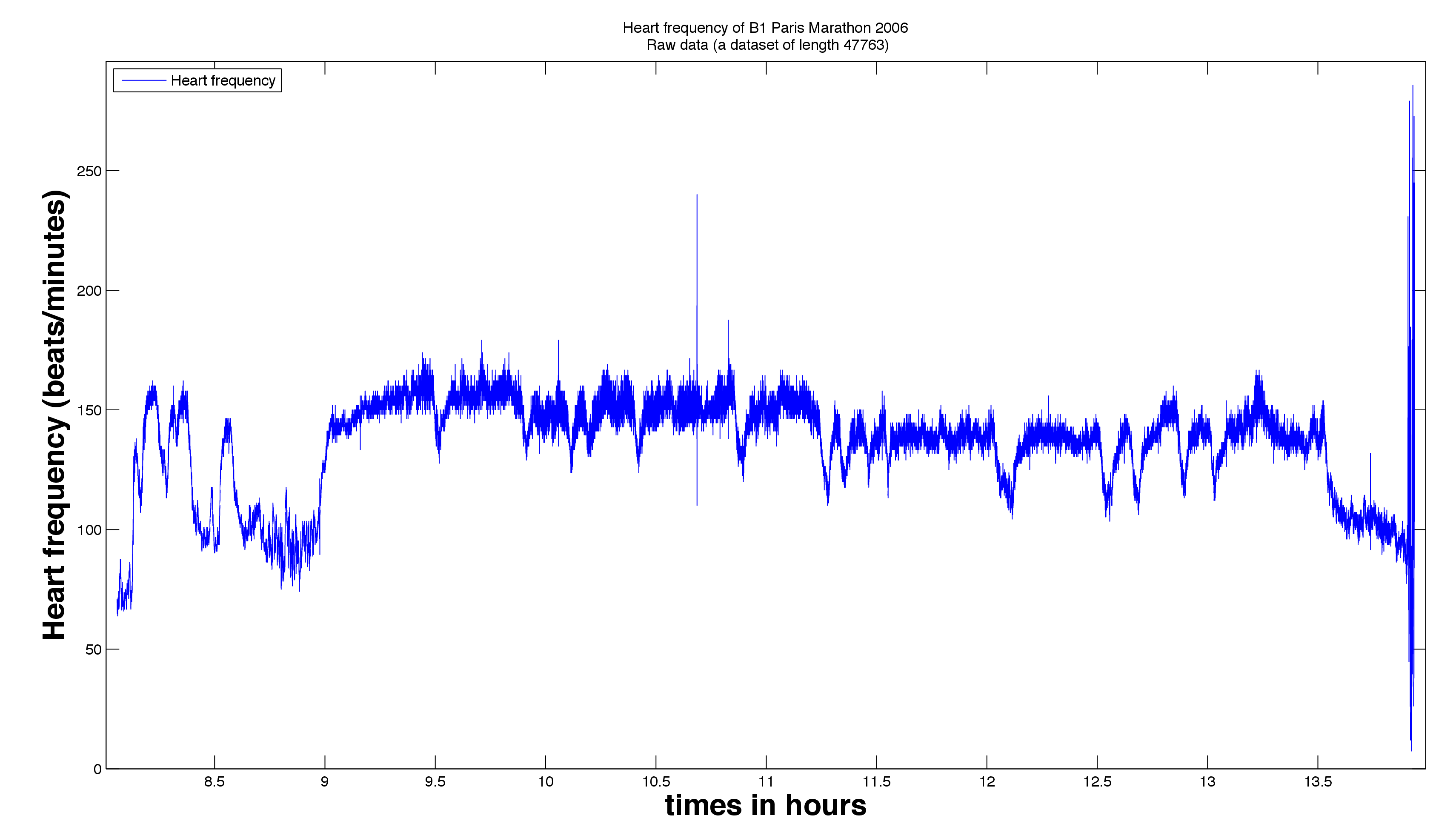}
\end{tabular}
\caption{\emph{Raw heart rate  time series of a marathon runner}}\label{fig14}
\end{center}
\end{figure}
In Figure~\ref{fig14}, Figure~\ref{fig15} and Figure~\ref{fig16}, we give two examples: 
In Figure~\ref{fig14}, the hearbeat time series correspond to  a marathon runner. This raw dataset  has been recorded by  
the team UBIAE (U. 902, INSERM and \'Evry G\'enopole) during Paris Marathon 2006.
This work is part of the project {\sc "Physio\-stat"} (2009-2011)  which is
supported by device grants from Digitéo and Région Île-de- France.
 The  hearbeat time series used in   Figure~\ref{fig16} 
     corresponds to a  day in the life of a shift worker. It has been kindly given by Gil Boudet and Professor Alain Chamoux (Occupational Safety service of Clermont Hospital). 
     
 Data are then been preprocessed by using the  {\sc "tachogram cleaning" } algorithm developed by Nadia Khalfa  and P.R.~Bertrand at INRIA Saclay in 2009.
{\sc "Tachogram cleaning" } cancelled aberrant data, based on physiological considerations rather than statistical procedure.
Next, segmentations of the cleaned heart beat time series can be  obtained 
by using the software "{\it In Vivo} Tachogram Analysis ({\it InVi}TA)" developed by P.R.~Bertrand at INRIA Saclay in 2010, see Figure~\ref{fig15} and  Figure ~\ref{fig16} below. 
 \begin{figure}[htbp]
\begin{center}
\begin{tabular}{c}
\includegraphics[width=15cm,height=6cm]{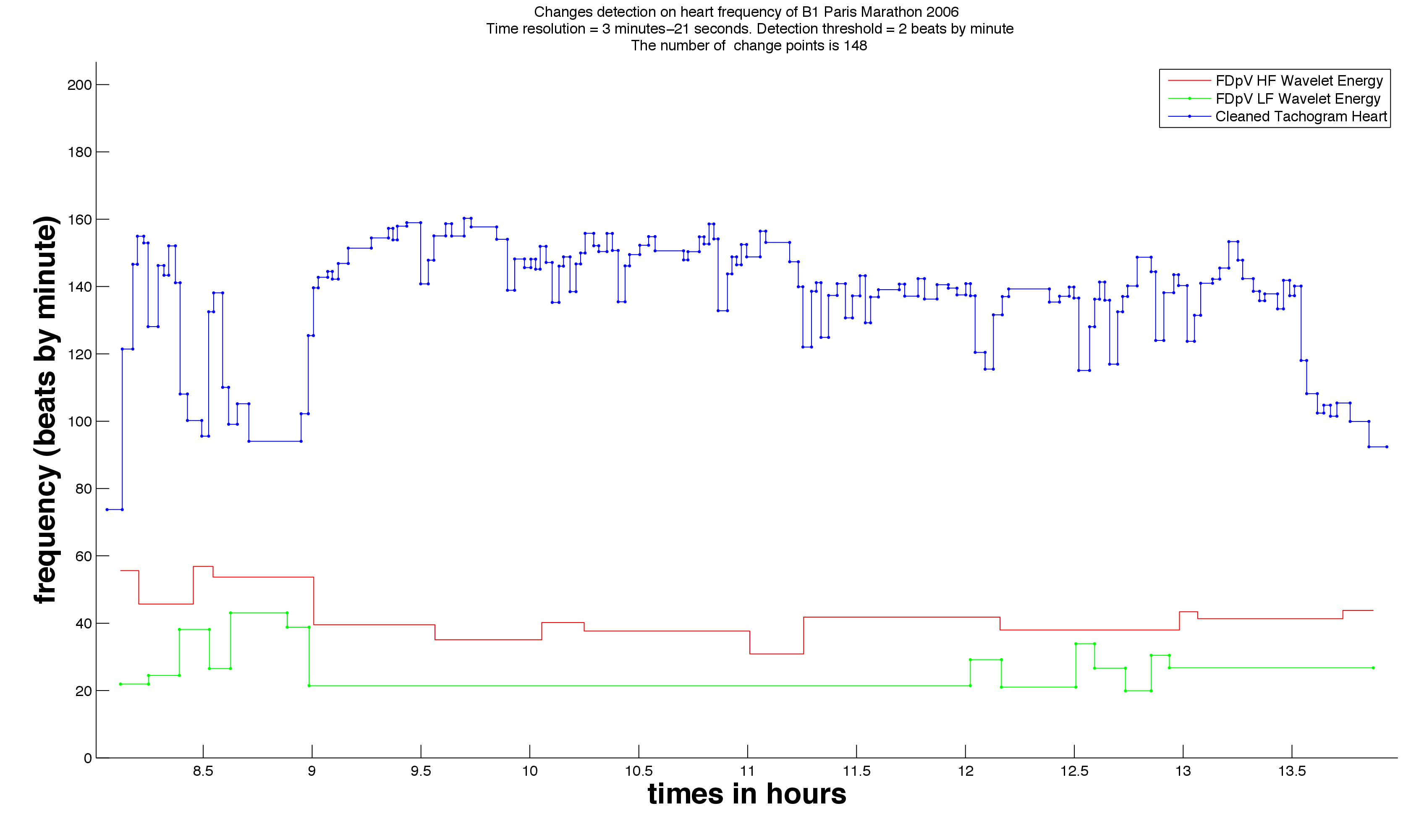}
\end{tabular}
\caption{\emph{Segmentation of heartbeat time series of a marathon runner. Blue: Heart rate  segmented by FDpV. 
Red: Low Frequency wavelet energy segmented by FDpV. Green: High Frequency wavelet energy segmented by FDpV
}}\label{fig15}
\end{center}
\end{figure}
Let us comment Fig.~\ref{fig15}  and Fig.~\ref{fig16}: first, for readability, we give the times in hours and minutes in abscissa and  the heart rate rather than RR-interval  in ordinate. We recall the equation: 
$\quad  Heart Rate = 60/RR~Interval\quad$ and we stress that all the computations are done on RR-interval. Secondly, in Figure~\ref{fig15}, we notice beginning and end of the marathon race, but also training before the race, and small breaks during the race. 
The same FDpV compression technology has been applied to wavelet energy corresponding to High Frequency, resp. Low Frequency band. Energy into these two frequencies bands is interpreted by cardiologists as corresponding to heart rate regulation by sympathetic and ortho-sympathetic systems. We refer to \citet{Ayache:Bertrand:2011} and 
\citet{Khalfa:etal:2011} for detailed explanations. In this paper, we just want to show that FDpV detects changes on HF and LF energy at the beginning of the race.
\begin{figure}[htbp]
\begin{center}
\begin{tabular}{c}
\includegraphics[width=16cm,height=6cm]{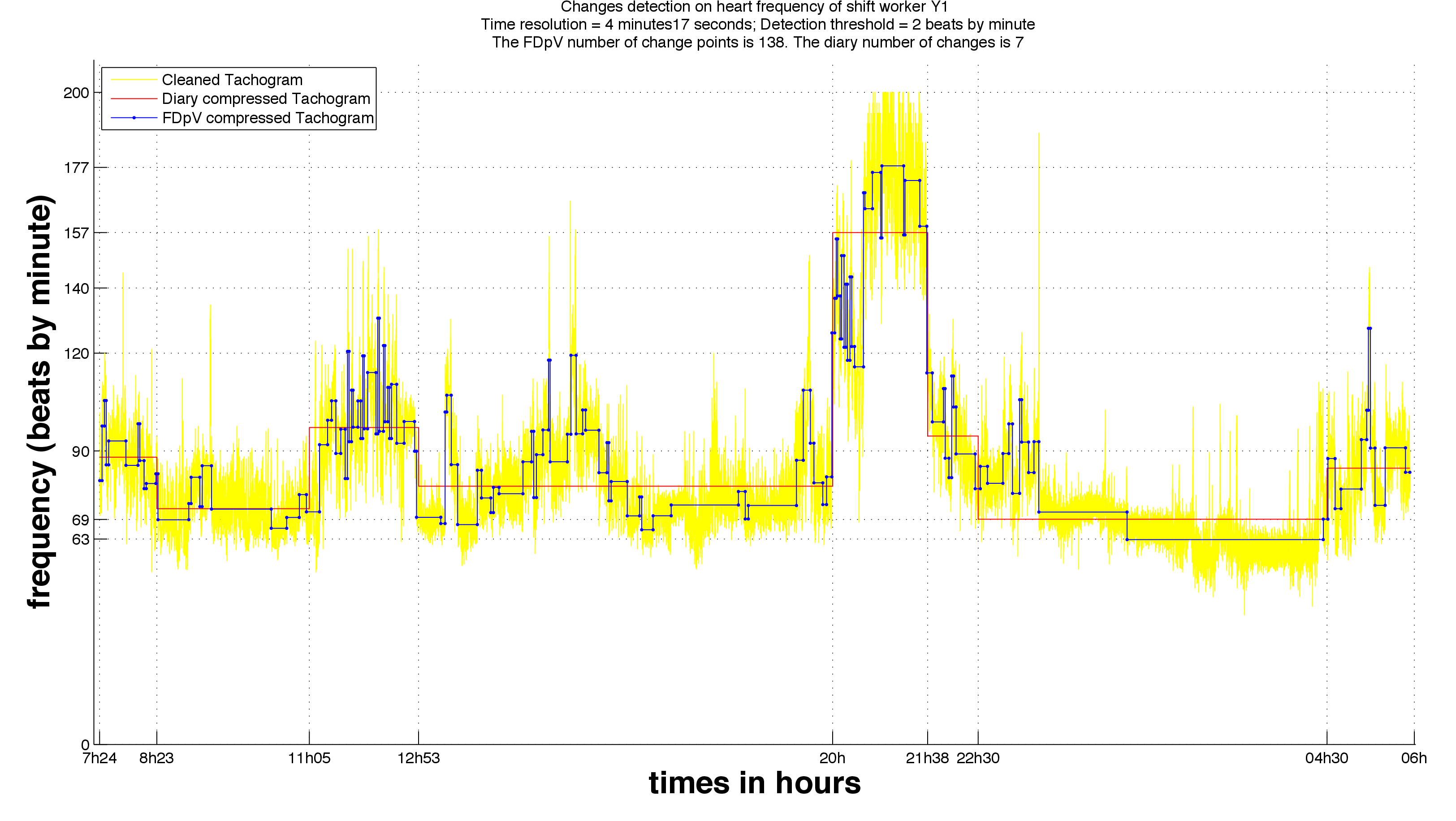}
\end{tabular}
\caption{\emph{Segmentation of heartbeat time series of shift worker Y1. Yellow: cleaned heart rate. Blue: heart rate segmented by FDpV. Red: Heart rate segmented following the dairy.}}\label{fig16}
\end{center}
\end{figure}
Fig.~\ref{fig16} is concerned with  heart beats of a shift worker Y1 through out a day its the life.
The shift worker Y1 has manually reported changes of activity on a diary, as shown in Table \ref{tab}.
\begin{table}[h!]
\begin{center}
\begin{tabular}{|l||c|c|c|c|c|c|c|}
\hline
 Time & 7h24-8h23 & 8h23-11h05 & 11h05-12h53 & 12h53-20h & 20h-21h38 & 22h30-4h30 \\
\hline
Activity & Task 1 & Task 2 & Picking & Free afternoon & {\bf playing football} & Sleeping\\
\hline
\end{tabular}
\caption{A shift worker's activities record}
\label{tab}
\end{center}
\end{table}
In yellow, we have plotted the cleaned heart rate time series. In red, we have plotted the segmentation resulting from manually recorded diary. In blue, we have plotted the automatic segmentation resulting from FDpV method.
Note that  the computation time is 15 seconds for 120,000 data with a code written in Matlab in a 2.8 GHz processor. Moreover, the FDpV segmentation is more accurate than the manual one.
For instance, Y1 has reported {\em "Football training"} from 20h to 21h38, which is supported by our analysis. But with FDpV method, we can see more details of the training namely the warming-up, the time for coach's recommendations,  and the football game with two small breaks.

Let us remark once again that both heartbeat time series and wavelet coefficient series do not fulfill the assumption of independency. However, FDpV segmentation provides accurate information on the heart rate, letting suppose a certain robustness of the method with respect to the stochastic model. 
 
In this case, the  FD-pV method has the advantage of being a fast and automatic procedure of segmentation of a large dataset, on the mean in this example, but possibly on the slope or the variance.  Combined with the  {\sc "tachogram cleaning" } algorithm, we then have an entirely automatic procedure to obtain apparently homogeneous segment. The next step 
will be to detect  change on hidden structural parameters. First attempts in this direction are exposed in \citet{Ayache:Bertrand:2011} and \citet{Khalfa:etal:2011}, but should be completed by 
forthcoming studies.

 \subsection*{Changes in the average daily volume}

\begin{figure}[htbp]
\begin{center}
\begin{tabular}{c}
\includegraphics[width=15cm,height=7cm]{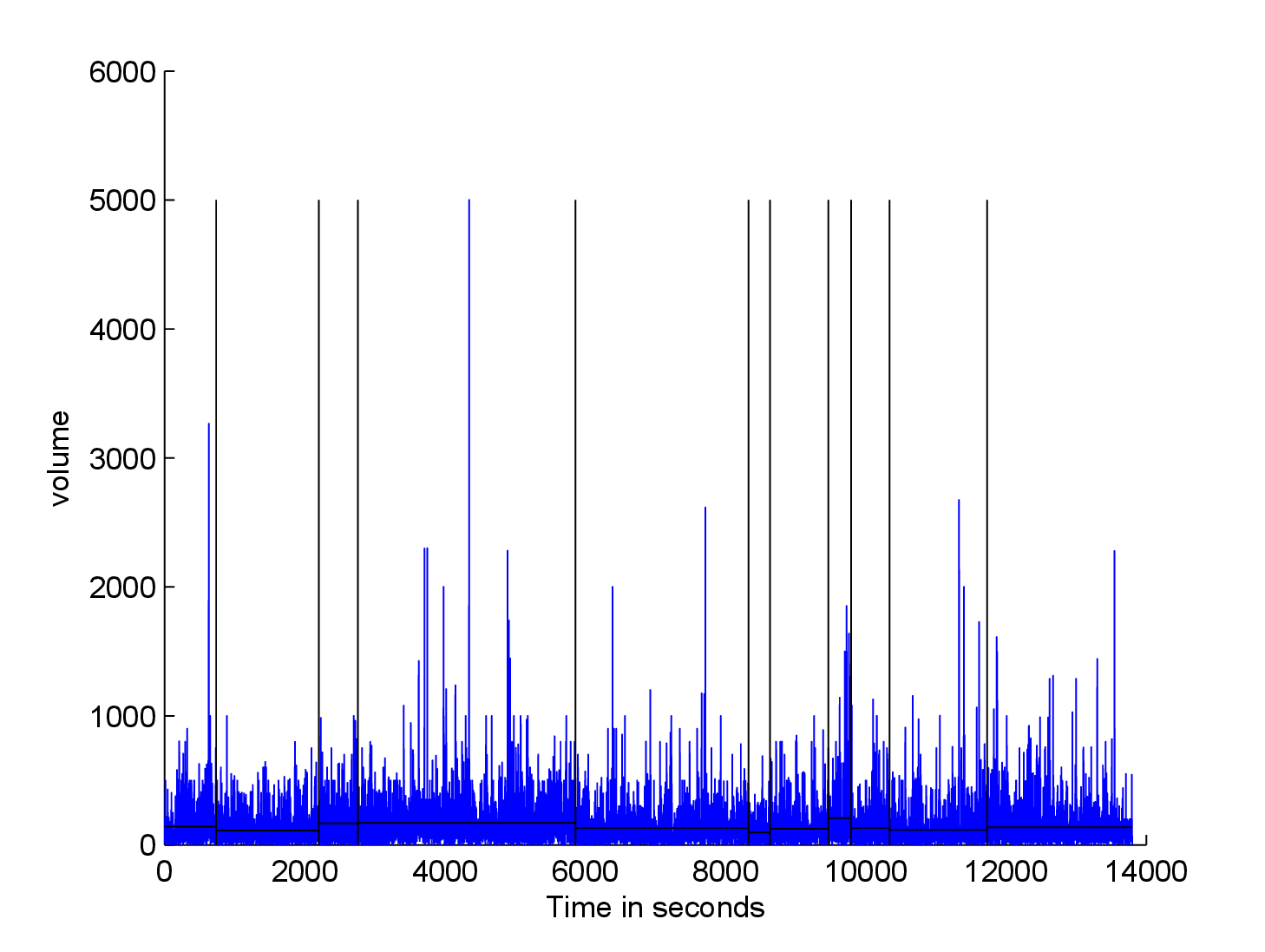}
\end{tabular}
\caption{\emph{Segmentation of the daily volume of Carbone
Lorraine compagny observed during 02 January 2009.}}\label{fig17}
\end{center}
\end{figure}

Trading volume represents number of shares or contracts traded in a financial market during a specific period. Average traded volume is an important indicator in technical analysis as it is used to measure the worth of a market move. If the markets move significantly up or down, the perceived strength of that move depends on the volume of trading in that period. The higher the volume during that price move, the more significant the is move.
Therefore, the detection of abrupt changes in the average
daily volume provide relevant information for financial engineer, trader, etc. Then, we consider here a daily volume of \emph{Carbone Lorraine} compagny observed during 02 January 2009. These data have been kindly given by Charles-Albert Lehalle from Cr\'edit Agricole Cheuvreux, Groupe CALYON (Paris). The results obtained with our algorithm for $A=300$, $p_1^{*}=0.05$ and
$p_2^{*}=10^{-5}$ are illustrated in Figure~\ref{fig17}. It
appears that FDp-V procedure detects major changes observed after each huge variations. In future works, we will investigate sequential detection of change points in the average daily volume in connection with the worth of a market move.

\begin{Remark}
We note that in our procedure, the observations are assumed to be independent. However, by Durbin and Watson test, we show that this condition is not checked, and we obtain a global correlation $\rho_{\text{Heartbeat}}=0.9982$ for the heartbeat data, and  $\rho_{\text{volume}}=0.4663$ for the financial time series. But this is not a restriction, on the contrary, this shows the robustness of our algorithm. In future works, this algorithm may be adapted to dependent data.
\end{Remark}

\section*{CONCLUSION}
It appears that both methods, namely FDp-V and PLSC, give
right results with practically the same precision. But, when we
compare the complexity, we remark that the FDp-V method is less
expensive in terms of time and memory complexity. Consequently,
 FDp-V method is faster (time) and cheaper (memory), and so it is
more adapted to segment random signals with large or huge datasets.

In future works, we will develop the Filtered Derivative with
p-value method in order to detect abrupt changes in parameters of
weakly or strongly dependent time series. In particularly, we will
consider the detection problem on the Hurst parameter of
multifractional Brownian motion
and apply it to physiological data as in \citet{Billat:Hamard:Meyer:Westfreid:2009}.
Let us also mention that  the
FDp-V method is based on sliding window and could be adapted to
sequential detection, see for instance \citet{Bertrand:Fleury:2008,Bertrand:Fhima:2009}.

\section{PROOF}\label{sec6}
\emph{\textbf{\large{Proof of Theorem~\ref{th1}.}}} Under the null
hypothesis $(H_0)$ the filtered derivative can be
rewritten as follows\\
$$\frac{\sqrt{A_n}}{\sigma}D_1(k,A_n)=\frac{S_{k-A_n}-2S_k+S_{k+A_n}}{\sqrt{A_n}},$$
where $\displaystyle S_k=\sum_{j=1}^{k}\xi_j$, with $(\xi_j)_{j
\in [1,n]}$ a sequence of i.i.d r.v such as $\E[\xi_1]=0$,
$\E[\xi_1^2]=1$, $S_0=0$.

\text{}\\
To achieve our goal, we state three lemmas. First, we show in
Lemma~\ref{lem:Lemme 1} that if a positive sequence $(a_n)$ has
the following asymptotic distribution
$$\lim\limits_{\substack{n \to +\infty}} \PP \left ( a_n \leq c_n(x) \right )=\exp(-2e^{-x}),$$
where $c_n(.)$ is defined by $\eqref{cnx}$ and if there is a
second positive sequence $(b_n)$ which converges almost surely
(a.s) to $(a_n)$  with rate of convergence of order $\mathcal{O
}\left( \sqrt{\log n} \right)$, then $(b_n)$ has the same
asymptotic distribution as $(a_n)$. We then prove in
Lemma~\ref{lem:Lemme 2} that under one of the assumptions
$(\mathcal{A}_i)$ with $i \in \{1,2,3\}$, the maximum of the
increment $\displaystyle \frac{\sqrt{A_n}}{\sigma} \left |
D_1(k,A_n) \right |$ converges a.s. to the maximum of discrete
Wiener process' increment with rate $\mathcal{O }\left( \sqrt{\log
n} \right)$. We show in Lemma~\ref{lem:Lemme 3} that the
maximum of \textbf{discrete} Wiener process' increment converges
a.s. to the maximum of \textbf{continuous} Wiener process'
increment with rate $\mathcal{O }\left( \sqrt{\log n} \right)$.
Then, by applying \citet[Theorem 5.2, p.
594]{Qualls:Watanabe:1972}, we deduce the asymptotic distribution
of the maximum of continuous Wiener process' increment. Finally,
by combining these results, we get
directly $\eqref{lim:mu}$.\\
To begin with, let us state the first lemma.

\begin{lemm}
\label{lem:Lemme 1}
Let $(a_n)$ and $(b_n)$ two sequences of positive random variables and we denote $\eta_n=|a_n-b_n|$.\\
We assume that
\begin{enumerate}
    \item $\lim\limits_{\substack{n \to +\infty}} \PP \left ( a_n \leq c_n(x) \right )=\exp(-2e^{-x}).$
    \item $\lim\limits_{\substack{n \to +\infty}} \eta_n \sqrt{\log
    n} \stackrel{a.s}{=} 0.$
\end{enumerate}
where $c_n(.)$ is defined by $\eqref{cnx}$. Then
\begin{equation}
\lim\limits_{\substack{n \to +\infty}} \PP \left ( b_n \leq c_n(x)
\right )=\exp(-2e^{-x}). \label{lim:lem1}
\end{equation}
\begin{flushright}
$\lozenge$
\end{flushright}
\end{lemm}
\emph{\textbf{\large{Proof of Lemma~\ref{lem:Lemme 1}.}}} Without
any restriction, we can consider the case where $\eta_n > 0$. We
denote by $|\Delta x|$ an infinitesimally small change in $x$.
Then, for $n$ large enough and $|\Delta x|$ small enough, $c_n
\left ( x + |\Delta x| \right )$ and $c_n \left ( x - | \Delta x |
\right )$ satisfy the following inequalities
\begin{eqnarray}
c_n \left ( x + |\Delta x_n| \right ) & \ge & c_n(x) + \frac{|\Delta x_n|}{\sqrt{2 \log n}}, \label{ineg1} \\
c_n \left ( x - |\Delta x_n| \right ) & \le & c_n(x) -
\frac{|\Delta x_n|}{\sqrt{2 \log n}}. \label{ineg2}
\end{eqnarray}

\text{}\\
Following \citet{Chen:1988}, we supply a lower and an
upper bounds of $\displaystyle \PP \left ( b_n \leq c_n(x) \right
)$. On the one hand, by using the inequality $\eqref{ineg1}$, the
upper bound results from the following calculations
\begin{eqnarray*}
\displaystyle \PP \left ( b_n \leq c_n(x) \right ) &\le&
\displaystyle \PP \left ( b_n \leq c_n \left (x + |\Delta x|
\right) - \frac{|\Delta x|}{\sqrt{2 \log n}} \right ) \\
& \le & \displaystyle \PP \left ( a_n - \eta_n \leq c_n \left (x +
|\Delta x| \right) - \frac{|\Delta x|}{\sqrt{2 \log n}} \right )\\
& \le & \displaystyle \PP \left ( a_n \leq c_n \left (x +
|\Delta x| \right) \right) + \PP \left ( - \eta_n \leq - \frac{|\Delta x|}{\sqrt{2 \log n}} \right )\\
& = & \displaystyle \PP \left ( a_n \leq c_n \left (x + |\Delta
x| \right) \right) + \PP \left ( \eta_n \ge  \frac{|\Delta
x|}{\sqrt{2 \log n}} \right ).
\end{eqnarray*}
On the other hand, by using the inequality $\eqref{ineg2}$, the
lower bound results from analogous calculations
\begin{eqnarray*}
\displaystyle \PP \left ( a_n \leq c_n \left (x - |\Delta x_n|
\right) \right ) & \le &  \displaystyle \PP \left ( b_n \leq c_n (x) \right) + \PP
\left ( \eta_n \ge  \frac{|\Delta x_n|}{\sqrt{2 \log n}} \right ).
\end{eqnarray*}
Then, by putting together the two previous bounds, we obtain
\begin{eqnarray*}
\PP \left ( a_n \leq c_n \left (x - |\Delta x_n| \right) \right )
-  \PP \left ( \eta_n \ge  \frac{|\Delta x_n|}{\sqrt{2 \log n}}
\right )
&\leq& \PP \left ( b_n \leq c_n (x) \right)\\
&\leq& \PP \left ( a_n\leq c_n \left (x + |\Delta x_n| \right)
\right) + \PP \left ( \eta_n \ge \frac{|\Delta x_n|}{\sqrt{2 \log
n}} \right ).
\end{eqnarray*}
\text{}\\
Finally, by taking the limit in $n$ and as $|\Delta x|$ is arbitrary small, we deduce $\eqref{lim:lem1}$. This
finishes the proof of Lemma~\ref{lem:Lemme 1}. \hfill $\blacklozenge$
\text{}\\
\text{}\\
Now, we show that the maximum of $
\frac{\sqrt{A_n}}{\sigma} \left | D_1(k,A_n) \right |$ converges
a.s to the maximum of discrete Wiener process' increment with rate
of convergence of order $ \mathcal{O }\left( \sqrt{\log n}
\right)$. This is stated in Lemma~\ref{lem:Lemme 2}
below (which gives also corrections to previous results of \citet{Chen:1988}:
\begin{lemm}
\label{lem:Lemme 2} Let $(W_t,\ t \ge 0)$ be a standard
Wiener process and $ \displaystyle \left ( Z_{A_n}(q), 0 \le q \le
\frac{n}{A_n}-1 \right )$ be the discrete sequence defined by
\begin{equation} Z_{A_n}(q)  = \left\{
\begin{array}{ll}
    \displaystyle W_{q-1}-2W_q+W_{q+1} & \mbox{if } 1 \le q \le \frac{n}{A_n}-1,\\
    0 & \mbox{else }.
\end{array}
\right. \label{def:Z}
\end{equation}
Let $\displaystyle (X_i)_{i=1,\dots,n}$ be a sequence of
independent identically distributed random variables with mean
$\mu$, variance $\sigma^2$, $D_1$ be defined by $\eqref{def:D1}$,
and denote
$$
\displaystyle \eta_{1,n}= \left| \max_{ 0 \le k \le n-A_n}
\frac{\sqrt{A_n}}{\sigma} \left|D_1(k,A_n) \right | - \max_{ 0 \le
q \le \frac{n}{A_n}-1} \left|Z_{A_n}(q) \right | \right|.
$$
Moreover, we suppose that one of the assumptions
$(\mathcal{A}_i)$, with $i \in \{1,2,3\}$ is in force. Then there exists a Wiener process $(W_t)_t$ such that
\be\label{lim:eta1} \lim\limits_{\substack{n \to +\infty}}
\eta_{1,n} \sqrt{\log n} \stackrel{a.s}{=}0. \ee
\end{lemm}
\begin{flushright}
$\lozenge$
\end{flushright}
\text{}\\
\emph{\textbf{\large{Proof of Lemma~\ref{lem:Lemme 2}.}}}
We consider a new discrete sequence, $  \left (
B(k,A_n), 0 \le k \le n-A_n \right )$, obtained by scaling from
the sequence $  \left ( Z_{A_n}(q), 0 \le q \le
\frac{n}{A_n}-1 \right )$. It is defined as follows
$$    B(k,A_n)  = \left\{
\begin{array}{ll}
    \displaystyle \frac{W_{k-A_n}-2W_k+W_{k+A_n}}{\sqrt{A_n}} & \mbox{if } A_n \le k \le n-A_n\\
    0 & \mbox{else }
\end{array}
\right. \\
$$
Then $$ \eta_{1,n} \stackrel{\mathcal{D}}{=} \left|
\max_{ 0 \le k \le n-A_n} \frac{\sqrt{A_n}}{\sigma} \left|D(k,A_n)
\right | - \max_{ 0 \le k \le n-A_n} \left|B(k,A_n) \right |
\right|,$$ where the sign $\stackrel{\mathcal{D}}{=}$ means
equality in~law. Depending on which assumption $(\mathcal{A}_i)$
is in force, we have three different proofs:
\begin{enumerate}
    \item Assuming $\mathbf{(\mathcal{A}_1)}$.\\
    This is the simplest case. We can choose a standard Wiener process, $(W_t,\ t \ge 0)$, such that
     $S_k = W_k$ at all the integers  $k \in [0,n-A_n]$. Hence, $\frac{\sqrt{A_n}}{\sigma}
    D(k,A_n)=B(k,A_n)$. Then, we can deduce $\eqref{lim:eta1}$.
    \item
    Assuming $\mathbf{(\mathcal{A}_2)}$.\\
    We have
    $$\displaystyle 0 \leq \eta_{1,n} \leq \max_{ 0 \le k \le n-A_n} \left|D(k,A_n) -  B(k,A_n) \right
    | \leq \frac{4}{\sqrt{A_n}} \max_{ 0 \le k \le n-A_n} \left|S_k -  W_k \right
    |$$
    and after
    $$ 0 \leq \eta_{1,n} \sqrt{\log n} \leq \frac{4 (\log n)^{\frac 3 2}}{\sqrt{A_n}}
    \times \frac{\displaystyle \max_{ 0 \le k \le n-A_n} \left|S_k -  W_k \right |}{\log n}$$
    However, $ \lim\limits_{\substack{n \to +\infty}} \frac{(\log n)^3}{A_n}=0$ and
    according to \citet[Theorem 3, p.34]{Komlos:etal:1975}, there is a Wiener process, $(W_t,\ t \ge
0)$,
    such as $ \lim\limits_{\substack{n \to +\infty}}
    \frac{{\displaystyle \max_{ 0 \le k \le n}} \left|S_k -  W_k \right |}{\log n} < + \infty .$
    Then, we can deduce $\eqref{lim:eta1}$.
    \item
    Assuming $\mathbf{(\mathcal{A}_3)}$.\\
    By using the same tricks than in the proof of the case
    $(\mathcal{A}_2)$, we can show that
    $$\displaystyle 0 \leq \eta_{1,n} \sqrt{\log n} \leq \frac{4 n^{\frac 1 p}\sqrt{\log n}}{\sqrt{A}}
    \times \frac{\displaystyle \max_{ 0 \le k \le n-A_n} \left|S_k -  W_k \right |}{n^{\frac 1 p}}$$
    However, $ \displaystyle \lim\limits_{\substack{n \to +\infty}} \frac{n^{\frac{2}{p}}\log n}{A_n}=0$
    and according to \citet[Theorem 3, p.34]{Komlos:etal:1975}, there is a Brownian motion, $(W_t,\ t \ge 0)$,
    such as $\displaystyle \lim\limits_{\substack{n \to +\infty}}
    \frac{\displaystyle \max_{ 0 \le k \le n} \left|S_k -  W_k \right |}{n^{\frac 1 p}} < + \infty .$
    Then, we can deduce (\ref{lim:eta1}).
    \end{enumerate}
This finishes the proof of Lemma~\ref{lem:Lemme 2}. \hfill$\blacklozenge$
\text{}\\
In order to apply \citet[Theorem 5.2, p.
594]{Qualls:Watanabe:1972} theorem, we have to consider a
continuous version of the process $Z_{A_n}$. For this reason, we
define the continuous process $\displaystyle \left( Z_{A_n}(t), t
\in \left[ 0,\frac{n}{A_n}-1 \right] \right)$ such as
\be\label{procZA} Z_{A_n}(t)=W_{t-1}-2W_t+W_{t+1}. \ee Then, in
Lemma~\ref{lem:Lemme 3}, we show that the maximum of
$\displaystyle \left | Z_{A_n}(q) \right |$ converges a.s to the
maximum of $\displaystyle \left | Z_{A_n}(t) \right |$ with rate
of convergence of order $ \mathcal{O }\left( \sqrt{\log n}
\right)$.
\begin{lemm}
\label{lem:Lemme 3} Let $Z_{A_n}$ be defined by
\eqref{def:Z} and set $\displaystyle \eta_{2,n}= \left| \sup_{  t \in \left [
0, \frac{n}{A_n}-1 \right ]} \left|Z_{A_n}(t) \right | - \max_{ 0
\le q \le \frac{n}{A_n}-1} \left|Z_{A_n}(q) \right | \right|.$\\
\emph{Then} \be\label{lim:eta2} \lim\limits_{\substack{n \to
+\infty}} \eta_{2,n} \sqrt{\log n} \stackrel{a.s}{=}0. \ee
\begin{flushright}
$\lozenge$
\end{flushright}
\end{lemm}
\emph{\textbf{\large{Proof of Lemma~\ref{lem:Lemme 3}.}}} We have
$$ 0 \leq \eta_{2,n} \leq  4 \sup_{  t \in \left [
0, \frac{n}{A_n}-1 \right ]} \left|W_{t+1} -  W_t \right |
 \leq  4 \sup_{ t \in \left [ 0,
\frac{n}{A_n}-\frac{1}{A_n}\right ]} \sup_{  s \in \left [
0,\frac{1}{A_n} \right ]} \left|W_{t+s} -  W_t \right | .
$$This implies
$$0 \leq \eta_{2,n} \sqrt{\log n} \leq 4 \sqrt{\log n}
\sup_{ t \in \left [ 0, \frac{n}{A_n}-\frac{1}{A_n} \right ]}
\sup_{  s \in \left [ 0,\frac{1}{A_n} \right ]}
\left|W_{t+s} - W_t \right | $$
and after
$$ \PP \left (
\eta_{2,n} \sqrt{\log n} \ge \delta \right) \leq \PP \left (
\sup_{ t \in \left [ 0, \frac{n}{A_n}-\frac{1}{A_n} \right ]}
\sup_{ s \in \left [ 0,\frac{1}{A_n} \right ]}
\left|W_{t+s} - W_t \right | \geq \frac{\delta
\sqrt{A_n}}{4 \sqrt{\log n}} \frac{1}{\sqrt{A_n}} \right )$$
\text{}\\
According to \citet[Lemma~1.2.1, p.
29]{Csorgo:Revesz:1981}, and by taking $\displaystyle T=n/A_n$,
$\displaystyle h=1/A_n$, $\varepsilon=1$, $\displaystyle
\nu=\frac{\delta \sqrt{A_n}}{4 \sqrt{\log n}}$,
and $C$ a non negative real, we deduce that\\
$$\PP \left ( \eta_{2,n} \sqrt{\log n} \ge
\delta \right) \leq C \frac{n}{A_n} \exp \left ( - \frac{\delta^2
A_n}{48 \log n} \right ).$$
\text{}\\
Next, by using $\lim\limits_{\substack{n \to +\infty}}
\frac{(\log n)^2}{A_n}=0,$ we can deduce $$\sum_{n \ge 1}
\frac{n}{A_n}     \exp \left ( - \frac{\delta^2 A_n}{48 \log n}
\right ) < + \infty,$$ and after $$\sum_{n \ge 1} \PP \left (
\eta_{2,n} \sqrt{\log n} \ge \delta \right) < + \infty.$$ Then,
according to Borel-Cantelli lemma, we can deduce
$\eqref{lim:eta2}$. This finishes the proof of
Lemma~\ref{lem:Lemme 3}. \hfill$\blacklozenge$
\text{}\\
Next, we apply \citet[Theorem 5.2, p.
594]{Qualls:Watanabe:1972} to the continuous process
$\displaystyle \left | Z_{A_n}(t) \right |.$ We obtain the
asymptotic distribution of its maximum given by \be \label{lim:Zt}
\lim\limits_{\substack{n \to +\infty}} \PP \left ( \sup_{ t \in
\left [ 0, \frac{n}{A}-1 \right ]} \left|Z_{A_n}(t) \right | \leq
c_n(x) \right )=\exp(-2e^{-x}). \ee
\text{}\\
This result can be proved by applying Theorem 5.2 of Qualls and
Watanabe (1972) to the centered stationary
Gaussian process $\left (Z_{A_n}(t), t \ge 0 \right )$. The
covariance function of $Z_{A_n}(t)$ is given by
$$
\displaystyle \rho(\tau)= \frac{Cov \left (Z_{A_n}(t),Z_{A_n}(t+\tau)
\right )}{\sqrt{Var \left ( Z_{A_n}(t) \right)} \sqrt{Var \left
(Z_{A_n}(t+\tau) \right)}}=\left\{
\begin{array}{ll}
     1-|\tau| & \mbox{if } 0 \le |\tau| \le 1,\\
     \displaystyle -1+\frac{|\tau|}{2} & \mbox{if } 1 < |\tau| \le 2,\\
     0 & \mbox{if } 2 < |\tau| < +\infty. \\
\end{array}
\right.
$$
So, the conditions of Theorem 5.2 of \citet{Qualls:Watanabe:1972}  are
satisfied in the following way: $\alpha=1$, $H_{\alpha}=1$
(according to \citet[p. 77]{Pickands:1969}, and
$\widetilde{\sigma}^{-1}(x)=2x^{-2}$. This finishes the proof of
$\eqref{lim:Zt}$.

Eventually, by combining Lemma~\ref{lem:Lemme 1}, Lemma~\ref{lem:Lemme
3} and result $\eqref{lim:Zt}$, we get   the asymptotic distribution
of the maximum of the sequence $ \displaystyle \left (
|Z_{A_n}(q)|, 0 \le q \le \frac{n}{A_n}-1 \right )$. Then, by
using Lemma~\ref{lem:Lemme 1} and Lemma~\ref{lem:Lemme 2}, we
immediately get the distribution of the maximum of the filtered
derivative sequence $(|D_1(A_n,k)|, 0 \leq k \leq n-A_n ).$\hfill $\blacksquare$
\text{}\\
\textbf{\large{End of the proof of Theorem~\ref{th1}.}}\\
\text{}\\
\emph{\textbf{\large{Proof of Theorem~\ref{th2}.}}}
Fix $\varepsilon >0$, the key argument is to divide $\Omega$ into
two complementary events
\begin{equation*}
 \Omega_{1,n} = \left \{ |\sigma-\widehat{\sigma}_n| \log n \leq \sigma \varepsilon \right \}
 \qquad
 \mathrm{and}
 \qquad
 \Omega_{2,n} = \left \{ |\sigma-\widehat{\sigma}_n| \log n > \sigma \varepsilon \right \}.
\end{equation*}
Then, we have
\begin{eqnarray*}
\PP  \left( \max_{k \in [A_n:n-A_n]} |D_1(k,A_n)| \le
\frac{\widehat{\sigma}_n}{\sqrt{A_n}} c_n(x) \right) &=& \PP
\left( \max_{k \in [A_n:n-A_n]} |D_1(k,A_n)| \le
\frac{\widehat{\sigma}_n}{\sqrt{A_n}} c_n(x) \text{ and }
\Omega_{1,n} \right)
\\
&+& \PP  \left( \max_{k \in [A_n:n-A_n]} |D_1(k,A_n)| \le
\frac{\widehat{\sigma}_n}{\sqrt{A_n}} c_n(x) \text{ and }
\Omega_{2,n} \right)
\end{eqnarray*}
On the one hand, we remark that
$$\PP  \left( \max_{k \in [A_n:n-A_n]}
|D_1(k,A_n)| \le \frac{\widehat{\sigma}_n}{\sqrt{A_n}} c_n(x)
\text{ and } \Omega_{2,n} \right) \leq \PP  \left( \Omega_{2,n}
\right),$$ which combined with assumption $\eqref{lim:hyp1}$
implies that
 \be\label{lim:omega2}
 \lim\limits_{\substack{n
\to +\infty}} \PP  \left( \Omega_{2,n} \right) = 0. \ee On the
other hand, for all $\omega \in \Omega_{1,n}$, we have
$\displaystyle\widehat{\sigma}_n=\sigma (1+\lambda_n(\omega))$
with $\displaystyle |\lambda_n(\omega)| \le
\frac{\varepsilon}{\log n}.$ Therefore, $\displaystyle
\lim\limits_{\substack{n \to +\infty}} \lambda_n
\stackrel{a.s}{=}0.$ Next, by setting $\displaystyle a_n =
\frac{\sqrt{A_n}}{\sigma}\max_{k \in [A_n:n-A_n]} |D_1(k,A_n)|$,
we get \ban \PP  \left( \max_{k \in [A_n:n-A_n]} |D_1(k,A_n)| \le
\frac{\widehat{\sigma}_n}{\sqrt{A_n}} c_n(x) \text{ and }
\Omega_{1,n} \right) &=& \PP  \left( a_n \leq
c_n(x)(1+\lambda_n(\omega)) \right)
\\
&=&
 \PP  \left( a_n -\eta_{3,n} \leq c_n(x) \right)
\ean with $\displaystyle \eta_{3,n}=c_n(x)\lambda_n$. Therefore,
after having checked that \be\label{lim:eta3}
\lim\limits_{\substack{n \to +\infty}} \eta_{3,n} \sqrt{\log n}
\stackrel{a.s}{=}0, \ee we can apply Lemma~\ref{lem:Lemme 1} which
combined  with Theorem~\ref{th1} implies that
$$
\PP  \left( \max_{k \in [A_n:n-A_n]} |D_1(k,A_n)| \le
\frac{\widehat{\sigma}_n}{\sqrt{A_n}} c_n(x) \text{ and }
\Omega_{1,n} \right)=\exp(-2e^{-x})$$ Eventually, combined with
$\eqref{lim:omega2}$ this implies $\eqref{lim:mu2}$. To finish the
proof, it just remains to verify that $\eqref{lim:eta3}$ is
satisfied. Indeed,
$$ \eta_{3,n} \sqrt{\log n} \leq \varepsilon \frac{c_n(x)}{\sqrt{\log n}},$$
and after having replaced $c_n(x)$ by its expression
$\eqref{cnx}$, we can easily verify that
$$\lim\limits_{\substack{n \to +\infty}} \eta_{3,n} \sqrt{\log n}
\stackrel{a.s}{=}0.$$ This finishes the proof of
Theorem~\ref{th2}.\hfill$\blacklozenge$
\text{}\\
\text{}\\
\emph{\textbf{\large{Proof of Corollary~\ref{cor2}.}}} Let $D_2$
and $\widehat{D}_2$ be defined respectively by $\eqref{def:D2}$
and $\eqref{def:D2hat}$, and set
$$
\eta_{4,n}= \left| \max_{ 0 \le k \le n-A_n}
\frac{\sqrt{A_n}}{\nu} \left|D_2(k,A_n) \right | - \max_{ 0 \le k
\le n-A_n} \frac{\sqrt{A_n}}{\nu} \left|\widehat{D}_2(k,A_n)
\right| \right|.
$$
\text{}\\
The key argument is to prove that
\begin{equation}
\lim\limits_{\substack{n \to +\infty}} \eta_{4,n} \sqrt{\log n}
\stackrel{a.s}{=}0. \label{lim:eta4}
\end{equation}
by using assumption $\eqref{condMu}$. Then, we apply
Lemma~\ref{lem:Lemme 1} which combined with Corollary~\ref{cor1}
implies $\eqref{lim:var2}$. So, to finish the proof we must verify
$\eqref{lim:eta4}$.\\
We have
$$
0 \leq \eta_{4,n} \leq \frac{\sqrt{A_n}}{\nu}  \max_{ 0 \le k \le
n-A_n} \left|D_2(k,A_n) - \widehat{D}_2(k,A_n) \right |
$$
which implies
$$
0 \leq \eta_{4,n} \leq \frac{\sqrt{A_n}}{\nu} \max_{ 0 \le k \le
n-A_n}  \left | |\mu-\widehat{\mu}_k|^2 -
|\mu-\widehat{\mu}_{k-A_n}|^2 \right |
$$
and after
$$
0 \leq \eta_{4,n} \sqrt{\log n} \leq \frac{2}{\nu} \max_{0 \le k
\le n-A_n} |\mu-\widehat{\mu}_k|^2 \sqrt{A_n \log n}.
$$
\text{}\\
Therefore, by using condition $\eqref{condMu}$, we get
$$
\lim\limits_{\substack{n \to +\infty}} \eta_{4,n} \sqrt{\log n}
\stackrel{a.s}{=}0.
$$
This finishes the proof of Corollary~\ref{cor2}. \hfill $\blacksquare$
\text{}\\

\section*{PROOF FOR LINEAR REGRESSION}
\emph{\textbf{\large{Proof of Theorem~\ref{th3}.}}} First we note
that, under the null hypothesis $(H_0)$, the sequence
$(D_3(k,A_n))_{A_n \leq k \leq n-A_n}$ satisfy
$$A D_3(k,A_n)= S_k^{-2} \sum_{j=k+1}^{k+A_n} \left ( X_j - \overline{X}_k \right ) \varepsilon_j
-S_{k-A_n}^{-2} \sum_{j=k-A_n+1}^{k} \left ( X_j -
\overline{X}_{k-A_n} \right ) \varepsilon_j
$$
where
\ban
\overline{X}_k = A_n^{-1} \sum_{j=k+1}^{k+A_n} X_j &\text{
and } &S_k^{2}= A_n^{-1} \sum_{j=k+1}^{k+A_n} \left ( X_j -
\overline{X}_k \right )^2,
\ean are respectively the empirical mean
and the
empirical variance of $X$ on the (sliding) box $[k+1,\,k+A_n]$.
By using the definition $\eqref{absXi}$, we see that
\begin{eqnarray*}
\overline{X}_k = \Delta \left ( k + \frac{A_n+1}{2} \right )
&\quad\mathrm{and}\quad&
S_k^{2} = \Delta^2 \left ( \frac{A_n^2-1}{12} \right ).
\end{eqnarray*}
Therefore, we can deduce
$$D_3(k,A_n)=\frac{12}{\Delta \times A_n \times (A_n^2-1)} \sum_{j=k-A_n+1}^{k+A_n} \gamma(j-k,A_n)\varepsilon_j$$
where \begin{equation} \gamma(i,A_n)= \left\{
\begin{array}{ll}
    i-\frac{A_n+1}{2} & \mbox{if } i > 0,\\
    -i-\frac{A_n-1}{2} & \mbox{if } i \leq 0.\\
\end{array}
\right. \label{gamma}
\end{equation}
Remark that the mean and the variance of the Gaussian
sequence $(D_3(k,A_n))_{A_n \leq k \leq n-A_n}$ verify
\begin{eqnarray*}
\E \left [ D_3(k,A_n) \right ]= 0
&\quad\mathrm{and}\quad&
Var \left [ D_3(k,A_n) \right ] = \frac{24 \sigma^2}{\Delta^2
\times A_n \times (A_n^2-1)}.
\end{eqnarray*}
Moreover the variance does not depend on $k$. Then, $(D_3(k,A_n))_{A_n \leq k \leq n-A_n}$  is a centered
stationary Gaussian sequence. Next, theorem \ref{th3} becomes
an application of \citet[Theorem 2.1, p.
3]{Csaki:Gonchigdanzan:2002} which gives the asymptotic
distribution of the maximum of standardized stationary Gaussian
sequences $(Z_k)_{k \geq 1}$ with covariance $\ds r_n=\cov
\left(Z_1,Z_{n+1} \right)$ under the condition $\ds r_n
\log (n) \rightarrow 0$.
Let us define the standardized version of $D_3$  as
\begin{equation}\label{def:D3std}
D_3^{\text{Std}}(k,A)=\frac{D_3(k,A)}{\sqrt{\text{Var} \left [
D_3(k,A) \right ]}} \qquad \mathrm{for} \;A_n \leq k \leq n-A_n,
\end{equation}
and its covariance
\ban
 \Gamma (k_1,k_2)&=& \cov
\left( D_3^{\text{Std}}(k_1,A),D_3^{\text{Std}}(k_2,A) \right).
\ean
The following lemma provides the value of the covariance:
\begin{lemm}
\label{lem:Lemme 4} Let $\left(D_3^{\text{Std}}(k,A_n)\right)_{A_n
\leq k \leq n-A_n}$ the standardized stationary Gaussian sequence
defined by $\eqref{def:D3std}$. Then, its covariance matrix
denoted $\left(\Gamma(k_1,k_2) \right)_{1 \leq k_1,k_2 \leq n}$ is
given by
$$
\Gamma(k_1,k_2)=  \frac{6}{A_n (A_n^2-1)}
 \times \left\{
\begin{array}{ll}
    f_1 \left ( |k_2-k_1|,A_n \right ) & \mbox{if } 0 \le |k_2-k_1| < A_n\\
    f_2\left ( |k_2-k_1|,A_n \right )  & \mbox{if } A_n \le |k_2-k_1| \le 2A_n-1\\
    0              & \mbox{if } 2A_n-1 < |k_2-k_1| < +\infty \\
\end{array}
\right.
$$
where
$$  f_1 \left (p,A_n \right )= \displaystyle \frac{1}{6}
(A_n-p)(A_n^2-2A_{n}p-2p^2-1) + \frac{1}{12} p
(3A_n^2+2p^2-6A_{n}p+1)$$ and
$$f_2 \left (p,A_n \right )= \displaystyle - \frac{1}{12}
(2A_n-p)(A_n^2-10A_{n}p-2p^2-1).$$
\begin{flushright}
$\lozenge$ \end{flushright}
\end{lemm}
\text{}\\
\emph{\textbf{\large{Proof of Lemma~\ref{lem:Lemme 4}.}}} First we
note that, by symmetry property of covariance matrix, we can
restrict ourselves to the  case $k_2-k_1 \geq 0$. Next, let us
 distinguish three different expressions of the covariance
according to the value of $k_2-k_1$.
\begin{itemize}
    \item If $k_2-k_1 > 2A_n-1$, then $\Gamma(k_1,k_2)=0.$
    \item If $0 \leq k_2-k_1 \leq 2A_n-1$, then
    $$ C^{-1}\Gamma(k_1,k_2)= \sum_{j=k_2-A_n+1}^{k_1+A_n} \gamma(j-k_1,A_n)
    \gamma(j-k_2,A_n) \text{ with } C = \frac{6}{A_n (A_n^2-1)}$$
    By replacing $j-k_1$ with $u$, we get
    $$\displaystyle C^{-1}\Gamma(k_1,k_2)=\sum_{u=-A_n+1}^{A_n-k_2+k_1} \gamma(u,A_n)
    \gamma(u+k_2-k_1,A_n).$$
    But, in order to use formula $\eqref{gamma}$ 
    we must find the case where
    the sign of $u$ and $u+k_2-k_1$ remain constant where $u \in [-A+1,A-k_2+k_1]$.
    That is why we must distinguish the following  subcases
    \begin{itemize}
        \item If $0 \leq k_2-k_1 \leq A_n-1$, then
        \begin{eqnarray*}
        C^{-1}\Gamma(k_1,k_2) &=& \sum_{u=-A_n+1}^{k_1-k_2}
        \gamma(\underbrace{u}_{\leq 0},A_n)\gamma(\underbrace{u+k_2-k_1}_{\leq
        0},A_n)\\
         &+&\sum_{u=k_1-k_2+1}^{0} \gamma(\underbrace{u}_{\leq
         0},A_n)\gamma(\underbrace{u+k_2-k_1}_{>0},A_n)\\
         &+& \sum_{u=1}^{A_n-k_2+k_1} \gamma(\underbrace{u}_{>0},A_n)\gamma(\underbrace{u+k_2-k_1}_{>0},A_n)
        \end{eqnarray*}
        and then
        $C^{-1}\Gamma_{k_1,k_2}=f_1 \left ( k_2-k_1,A_n \right ).$
        \item If $A_n \leq k_2-k_1 \leq 2A_n-1$, then
        $$C^{-1}\Gamma(k_1,k_2)=\sum_{u=-A_n+1}^{A_n-k_2+k_1}
        \gamma(\underbrace{u}_{\leq 0},A_n) \gamma(\underbrace{u+k_2-k_1}_{>0},A_n)$$
        which implies
        $C^{-1} \Gamma(k_1,k_2)= f_2 \left ( k_2-k_1,A_n \right
        ).$
    \end{itemize}
\end{itemize}
This finishes the proof of Lemma~\ref{lem:Lemme 4}.\hfill$\blacklozenge$
\text{}\\

Hence, by applying  \citet[Theorem 2.1, p.
3]{Csaki:Gonchigdanzan:2002} to the
sequence $\left(D_3^{\text{Std}}(k,A_n)\right)_{A_n \leq k \leq
n-A_n}$, we can deduce that
\begin{equation}
\lim\limits_{\substack{n \to +\infty}} \PP \left( \max_{k \in
[A_n:n-A_n]} |D_3^{\text{Std}}(k,A_n)| \le d_n(x)  \right) = \exp(-2e^{-x}).
\end{equation}
Then, by using $\eqref{def:D3std}$, we obtain $\eqref{lim:a}$.
This finishes the proof of Theorem~\ref{th3}. \hfill$\blacksquare$
\text{}\\
\text{}\\
\emph{\textbf{\large{Proof of Corollary~\ref{cor3}.}}} First we
note that, under the null hypothesis $(H_0)$, the Filtered
Derivative applied to the intercept satisfy
$$
D_4(k,A_n)=A_n^{-1} \sum_{j=k+1}^{k+A_n} e_j - A_n^{-1}
\sum_{j=k-A_n+1}^{k} e_j.
$$
Therefore, $\left(D_4(k,A_n)\right)_{A_n \leq k \leq
n-A_n}$ corresponds to a sequence of Filtered Derivative of the mean, applied to the particular case of i.i.d centered Gaussian r.v with known variance $\sigma^2$. Then, by applying Theorem~\ref{th1} under assumption
$\mathcal{A}_1$, we obtain $\eqref{lim:b}$. This
finishes the proof of Corollary~\ref{cor3}. \hfill $\blacksquare$
\text{}\\
\text{}\\
\emph{\textbf{\large{Proof of Corollary~\ref{cor4}.}}} Let $D_4$
and $\widehat{D}_4$ be defined respectively by $\eqref{def:D4}$
and $\eqref{def:D4hat}$, and set
$$
\eta_{5,n}= \left| \max_{ 0 \le k \le n-A_n}
\frac{\sqrt{A_n}}{\sigma} \left|D_4(k,A_n) \right | - \max_{ 0 \le
k \le n-A_n} \frac{\sqrt{A_n}}{\sigma} \left|\widehat{D}_4(k,A_n)
\right| \right|.
$$
The key argument is to prove that
\begin{equation}
\lim\limits_{\substack{n \to +\infty}} \eta_{5,n} \sqrt{\log n}
\stackrel{a.s}{=}0. \label{lim:eta5}
\end{equation}
By using assumption $\eqref{condA}$. Then, we apply
Lemma~\ref{lem:Lemme 1} which combined with Corollary~\ref{cor3}
implies $\eqref{lim:b2}$. So, to finish the proof we must verify
$\eqref{lim:eta5}$.
Next, we have
$$
0 \leq \eta_{5,n} \leq \frac{\sqrt{A_n}}{\sigma}  \max_{ 0 \le k
\le n-A_n} \left|D_4(k,A_n) - \widehat{D}_4(k,A_n) \right |
$$
which implies
$$
0 \leq \eta_{5,n} \leq \frac{1}{\sqrt{A_n} \sigma}
|a-\widehat{a}_n| \max_{ 0 \le k \le n-A_n}  \left | \ds
\sum_{j=k+1}^{k+A_n} X_j - \sum_{j=k-A_n+1}^{k} X_j \right |.
$$
Then, by using $\eqref{absXi}$, we show that
$$
0 \leq \eta_{5,n} \leq \frac{\Delta_n}{\sqrt{A_n} \sigma}
|a-\widehat{a}_n| \max_{ 0 \le k \le n-A_n}
\underbrace{\sum_{j=k-A_n+1}^{k+A_n}j}_{=A_n^2}
$$
and after
$$
\eta_{5,n} \sqrt{\log n} \leq |a-\hat{a}_n| A_n^{\frac{3}{2}}
\Delta_n \sqrt{\log n}.
$$
Therefore, by using condition $\eqref{condA}$, we get
$$
\lim\limits_{\substack{n \to +\infty}} \eta_{5,n} \sqrt{\log n}
\stackrel{a.s}{=}0.
$$
This finishes the proof of Corollary~\ref{cor4} \hfill$\blacksquare$
\text{}\\
\text{}\\
\Acknowledgments{We are grateful to the  editors and the referees  for their helpful comments. A first version of this work was presented during the 
"International Workshop in Sequential Methodologies"  at UTT (Troyes, France, June  15-17, 2009). 
We would like to thank both organizers and  participants of this conference for  stimulating discussions and  suggestions.}

\renewcommand{\refname}{REFERENCES}

\bibliographystyle{apa}
\bibliography{OffLineChangePoint}

\begin{thebibliography}{}

\bibitem[\protect\astroncite{Antoch and
  Hu{\v{s}}kov{\'a}}{1994}]{Antoch:Huskova:1994}
Antoch, J. and Hu{\v{s}}kov{\'a}, M. (1994).
\newblock Procedures for the detection of multiple changes in series of
  independent observations.
\newblock In {\em Asymptotic statistics ({P}rague, 1993)}, Contrib. Statist.,
  pages 3--20. Physica, Heidelberg.

\bibitem[\protect\astroncite{Ayache and Bertrand}{2011}]{Ayache:Bertrand:2011}
Ayache, A. and Bertrand, P.~R. (2011).
\newblock Discretization error of wavelet coefficient for fractal like process.
\newblock {\em Advances in Pure and Applied Mathematics}, to appear.

\bibitem[\protect\astroncite{Bai and Perron}{1998}]{Bai:Perron:1998}
Bai, J. and Perron, P. (1998).
\newblock Estimating and testing linear models with multiple structural
  changes.
\newblock {\em Econometrica}, 66(1):47--78.

\bibitem[\protect\astroncite{Basseville and
  Nikiforov}{1993}]{Basseville:Nikiforov:1993}
Basseville, M. and Nikiforov, I.~V. (1993).
\newblock {\em Detection of abrupt changes: theory and application}.
\newblock Prentice Hall Information and System Sciences Series. Prentice Hall
  Inc., Englewood Cliffs, NJ.

\bibitem[\protect\astroncite{Ben-Yaacov and Eldar}{2008}]{BenYaacov:Eldar:2008}
Ben-Yaacov, E. and Eldar, Y.~C. (2008).
\newblock A fast and flexible method for the segmentation of acgh data.
\newblock {\em Bioinformatics}, 24:139--145.

\bibitem[\protect\astroncite{Benveniste and
  Basseville}{1984}]{Benveniste:Basseville:1984}
Benveniste, A. and Basseville, M. (1984).
\newblock Detection of abrupt changes in signals and dynamical systems: some
  statistical aspects.
\newblock In {\em Analysis and optimization of systems, {P}art 1 ({N}ice,
  1984)}, volume~62 of {\em Lecture Notes in Control and Inform. Sci.}, pages
  145--155. Springer, Berlin.

\bibitem[\protect\astroncite{Bertrand}{2000}]{Bertrand:2000}
Bertrand, P.~R. (2000).
\newblock A local method for estimating change points: the ``hat-function''.
\newblock {\em Statistics}, 34(3):215--235.

\bibitem[\protect\astroncite{Bertrand and Fhima}{2009}]{Bertrand:Fhima:2009}
Bertrand, P.~R. and Fhima, M. (2009).
\newblock Filtered derivative with p-value method for multiple change-points
  detection.
\newblock In {\em Proceeding of the 2nd International Workshop in Sequential
  Methodologies}.

\bibitem[\protect\astroncite{Bertrand and Fleury}{2008}]{Bertrand:Fleury:2008}
Bertrand, P.~R. and Fleury, G. (2008).
\newblock Detecting small shift on the mean by finite moving average.
\newblock {\em International Journal of Statistics and Management System},
  3:56--73.

\bibitem[\protect\astroncite{Billat
  et~al.}{2009}]{Billat:Hamard:Meyer:Westfreid:2009}
Billat, V.~L., Hamard, L., Meyer, Y., and Wesfreid, E. (2009).
\newblock Detection of changes in the fractal scaling of heart rate and speed
  in a marathon race.
\newblock {\em Physica A}, 388:3798--3808.

\bibitem[\protect\astroncite{Birg{\'e} and Massart}{2007}]{Birge:Massart:2007}
Birg{\'e}, L. and Massart, P. (2007).
\newblock Minimal penalties for {G}aussian model selection.
\newblock {\em Probab. Theory Related Fields}, 138(1-2):33--73.

\bibitem[\protect\astroncite{Brodsky and
  Darkhovsky}{1993}]{Brodsky:Darkhovsky:1993}
Brodsky, B.~E. and Darkhovsky, B.~S. (1993).
\newblock {\em Nonparametric methods in change-point problems}, volume 243 of
  {\em Mathematics and its Applications}.
\newblock Kluwer Academic Publishers Group, Dordrecht.

\bibitem[\protect\astroncite{Chen}{1988}]{Chen:1988}
Chen, X. (1988).
\newblock Inference in a simple change-point model.
\newblock {\em Scienta Sinica}, A31:654--667.

\bibitem[\protect\astroncite{Chopin}{2007}]{Chopin:2007}
Chopin, N. (2007).
\newblock Dynamic detection of change points in long time series.
\newblock {\em Annals of the Institute of Statistical Mathematics},
  59:349--366.

\bibitem[\protect\astroncite{Cs\'{a}ki and
  Gonchigdanzan}{2002}]{Csaki:Gonchigdanzan:2002}
Cs\'{a}ki, E. and Gonchigdanzan, K. (2002).
\newblock Almost sure limit theorems for the maximum of stationary gaussian
  sequences.
\newblock {\em Stat. Prob. Letters}, 58:195--203.

\bibitem[\protect\astroncite{Cs\"{o}rgo and
  Horv\'ath}{1997}]{Csorgo:Horvath:1997}
Cs\"{o}rgo, M. and Horv\'ath, L. (1997).
\newblock {\em Limit Theorem in Change-Point Analysis}.
\newblock J. Wiley, New York.

\bibitem[\protect\astroncite{Cs\"{o}rg\"{o} and
  R\'ev\'esz}{1981}]{Csorgo:Revesz:1981}
Cs\"{o}rg\"{o}, M. and R\'ev\'esz, P. (1981).
\newblock {\em Strong Approximations in Probability and Statistics}.
\newblock Akad\'{e}miai Kiad\"{o}, Budapest.

\bibitem[\protect\astroncite{Fearnhead and Liu}{2007}]{Fearnhead:Liu:2007}
Fearnhead, P. and Liu, Z. (2007).
\newblock On-line inference for multiple change points problems.
\newblock {\em Journal of the Royal Statistical Society}, 69:589--605.

\bibitem[\protect\astroncite{Fliess et~al.}{2010}]{Fliess:etal:2010}
Fliess, M., Join, C., and Mboup, M. (2010).
\newblock Algebraic change-point detection.
\newblock {\em Applicable Algebra in Engineering, Communication and Computing},
  21:131--143.

\bibitem[\protect\astroncite{Gombay and Serban}{2009}]{Gombay:Serban:2009}
Gombay, E. and Serban, D. (2009).
\newblock Monitoring parameter change in ar(p) time series models.
\newblock {\em Journal of Multivariate Analysis}, 100:715--725.

\bibitem[\protect\astroncite{Gu\'edon}{2007}]{Guedon:2007}
Gu\'edon, Y. (2007).
\newblock Exploring the state sequence space for hidden markov and semi markov
  chains.
\newblock {\em Computational Statistics and Data Analysis}, 51:2379--2409.

\bibitem[\protect\astroncite{Hu\v{s}kov\'a and
  Meintanis}{2006}]{Huskova:Meintanis:2006}
Hu\v{s}kov\'a, M. and Meintanis, S.~G. (2006).
\newblock Change point analysis based on the empirical characteristic functions
  of ranks.
\newblock {\em Sequential Analysis}, 25:421--436.

\bibitem[\protect\astroncite{Khalfa et~al.}{2011}]{Khalfa:etal:2011}
Khalfa, N., Bertrand, P.~R., Boudet, G., Chamoux, A., and Billat, V. (2011).
\newblock Heart rate regulation processed through wavelet analysis and change
  detection. some case studies.
\newblock {\em Submitted}.

\bibitem[\protect\astroncite{Kirch}{2008}]{Kirch:2008}
Kirch, C. (2008).
\newblock Bootstrapping sequential change-point tests.
\newblock {\em Sequential Analysis}, 27:330--349.

\bibitem[\protect\astroncite{Koml\'{o}s et~al.}{1975}]{Komlos:etal:1975}
Koml\'{o}s, J., Major, P., and Tusn\'{a}dy, G. (1975).
\newblock An approximation of partial sums of independent ${\rm rv}$'s and the
  sample ${\rm df}$.
\newblock {\em I. Z. Wahrscheinlichkeitstheorie und Verw. Gebiete},
  32:111--131.

\bibitem[\protect\astroncite{Lavielle and
  Moulines}{2000}]{Lavielle:Moulines:2000}
Lavielle, M. and Moulines, E. (2000).
\newblock Least-squares estimation of an unknown number of shifts in a time
  series.
\newblock {\em J. Time Ser. Anal.}, 21(1):33--59.

\bibitem[\protect\astroncite{Lavielle and
  Teyssi\`ere}{2006}]{Lavielle:Teyssiere:2006}
Lavielle, M. and Teyssi\`ere, G. (2006).
\newblock Detection of multiple change points in multivariate time series.
\newblock {\em Lithuanian Math. J.}, 46:287--306.

\bibitem[\protect\astroncite{Lebarbier}{2005}]{Lebarbier:2005}
Lebarbier, E. (2005).
\newblock Detecting multiple change-points in the mean of gaussian process by
  model selection.
\newblock {\em Signal Processing}, 85:717--736.

\bibitem[\protect\astroncite{Montgomery}{1997}]{Montgomery:1997}
Montgomery, D. (1997).
\newblock {\em Introduction to Statistical Quality Control, 3rd Edition.}
\newblock John Wiley \& Sons, New York.

\bibitem[\protect\astroncite{Pickands}{1969}]{Pickands:1969}
Pickands, J. (1969).
\newblock Asymptotic properties of the maximum in a stationary gaussian
  process.
\newblock {\em Trans. Amer. Math. Soc}, 145:75--86.

\bibitem[\protect\astroncite{Qualls and Watanabe}{1972}]{Qualls:Watanabe:1972}
Qualls, C. and Watanabe, H. (1972).
\newblock Asymptotic properties of gaussian processes.
\newblock {\em Ann. Math. Statist}, 43:580--596.

\bibitem[\protect\astroncite{R\'ev\'esz}{1990}]{Revesz:1990}
R\'ev\'esz, P. (1990).
\newblock {\em Random Walks in Random and non-random enviroments}.
\newblock World Scientific Publishing Company.

\bibitem[\protect\astroncite{Steinebach and
  Eastwood}{1995}]{Steinebach:Eastwood:1995}
Steinebach, J. and Eastwood, R. (1995).
\newblock On extreme value asymptotics for increments of renewal processes.
\newblock {\em Journal of Statistical Planning and Inference}, 45:301--312.

\bibitem[\protect\astroncite{{Task force of the European Soc. Cardiology and
  the North American Society of Pacing and
  Electrophysiology}}{1996}]{Task:Force:1996}
{Task force of the European Soc. Cardiology and the North American Society of
  Pacing and Electrophysiology} (1996).
\newblock {Heart rate variability. Standards of measurement, physiological
  interpretation, and clinical use.}
\newblock {\em Circulation}, 93:1043--1065.

\end{thebibliography}

\end{document}